\documentclass{article}

\usepackage[
backend=biber,
style=numeric-comp,
sorting=nty
]{biblatex}
\usepackage[margin=1in]{geometry}
\usepackage[utf8]{inputenc}
\usepackage{csquotes}
\usepackage{amsfonts,amsthm,amsmath,amssymb,accents,mathtools,physics} 
\usepackage{mathrsfs,amsbsy} 
\usepackage{enumitem}
\usepackage{subfig}
\usepackage{graphicx}
\usepackage{float}
\usepackage{xcolor}
\usepackage{algorithm}
\usepackage[noend]{algpseudocode}
\usepackage{hyperref}
\usepackage{authblk} 
\usepackage{textcomp} 
\usepackage{algorithm}
\usepackage{algorithmicx}
\usepackage{relsize}

\usepackage{bbm}
\usepackage{dsfont}

\numberwithin{equation}{section}

\usepackage{yhmath}
 \usepackage{booktabs} 

\newtheorem{thm}{Theorem}[section]

\newtheorem{defn}[thm]{Definition}
\newtheorem{rem}[thm]{Remark}

\theoremstyle{definition}
\newtheorem{exa}[thm]{Example}

\newcommand{\bit}{\begin{itemize}}
\newcommand{\eit}{\end{itemize}}
\newcommand{\beq}{\begin{equation}}
\newcommand{\eeq}{\end{equation}}
\newcommand{\be}{\begin{eqnarray}}
\newcommand{\ee}{\end{eqnarray}}
\newcommand{\beno}{\begin{eqnarray*}}
\newcommand{\eeno}{\end{eqnarray*}}

\hypersetup{
    colorlinks=true,
    linkcolor=blue,
    citecolor=blue,
    filecolor=blue,      
    urlcolor=blue,
}

\title{A high-order Eulerian-Lagrangian Runge-Kutta finite volume (EL-RK-FV) method for scalar nonlinear conservation laws}
\author[1]{Jiajie Chen}
\author[2]{Joseph Nakao}
\author[3]{Jing-Mei Qiu}
\author[4]{Yang Yang}
\affil[1]{Department of Biostatistics, Epidemiology, and Informatics, University of Pennsylvania Perelman School of Medicine, Philadelphia, PA, USA}
\affil[2]{Department of Mathematics and Statistics, Swarthmore College, Swarthmore, PA, USA}
\affil[3]{Department of Mathematical Sciences, University of Delaware, Newark, DE, USA}
\affil[4]{Department of Mathematical Sciences, Michigan Technological University, Houghton, MI, USA}
\date{}

\begin{document}

\maketitle


\begin{abstract}
\noindent We present a class of high-order Eulerian-Lagrangian Runge-Kutta finite volume methods that can numerically solve Burgers' equation with shock formations, which could be extended to general scalar conservation laws. Eulerian-Lagrangian (EL) and semi-Lagrangian (SL) methods have recently seen increased development and have become a staple for allowing large time-stepping sizes. Yet, maintaining relatively large time-stepping sizes post shock formation remains quite challenging. Our proposed scheme integrates the partial differential equation on a space-time region partitioned by linear approximations to the characteristics determined by the Rankine-Hugoniot jump condition. We trace the characteristics forward in time and present a merging procedure for the mesh cells to handle intersecting characteristics due to shocks. Following this partitioning, we write the equation in a time-differential form and evolve with Runge-Kutta methods in a method-of-lines fashion. High-resolution methods such as ENO and WENO-AO schemes are used for spatial reconstruction. Extension to higher dimensions is done via dimensional splitting. Numerical experiments demonstrate our scheme's high-order accuracy and ability to sharply capture post-shock solutions with large time-stepping sizes.
\end{abstract}



\section{Introduction}
\label{sec:introduction}
We are interested in numerically solving scalar (non)linear hyperbolic conservation laws of the form
\begin{equation}
	\label{eq: hyperbolicconservationlaw}
	u_t + \nabla\cdot\mathbf{F}(u)=0,\quad \mathbf{x}\in\mathcal{D},\quad t>0,
\end{equation}
where $\mathbf{F}(u)$ is a strictly convex flux function. Such equations appear in many physical systems, with various applications in plasma physics, continuum physics, meterology, traffic flow, fluid dynamics, and large-scale supply chains in economics. However, numerically solving equation \eqref{eq: hyperbolicconservationlaw} comes with several challenges: the time-stepping size is usually constrained since explicit time integrators are often utilized due to the nonstiff convection term, and shocks may develop in the solution even with sufficiently smooth initial data. As a first step to handling shock formations in hyperbolic conservation laws \eqref{eq: hyperbolicconservationlaw}, we consider Burgers' equation with scalar flux function $f(u)=u^2/2$ in each dimension. In doing so, we lay a framework that could be extended to general scalar conservation laws.

Eulerian-Lagrangian (EL) and semi-Lagrangian (SL) schemes \cite{russell2002overview,xiu2001semi} have become popular methods for solving linear and nonlinear hyperbolic conservation laws. This recent growth in popularity is driven by their ability to combine the strengths of both Eulerian and Lagrangian frameworks. Eulerian methods feature fixed numerical grids that allow high spatial resolution schemes to be incorporated, but the fixed mesh leads to a time-stepping constraint from the Courant-Friedrichs-Lewy (CFL) condition. Lagrangian methods transport information along the characteristics, hence allowing the solution to be evolved free of the CFL condition. However, Lagrangian methods are known to have a significant amount of noise, and following the trajectories that drive the system could lead to a clustering of nodes, e.g., around sharp gradients, and requires a remeshing step. Broadly speaking, EL and SL methods work on a stationary numerical grid to allow high spatial resolution schemes to be used, but the solution is evolved along approximate characteristics to relax the CFL condition and allow very large time-stepping sizes prior to the characteristics crossing each other, i.e., when shocks form. EL and SL methods have been developed in several frameworks over the past few decades: discontinuous Galerkin (DG) \cite{cai2017high, cai2018high, ding2020semi,hong2024conservative,qiu2011positivity,rossmanith2011}, finite difference (FD) \cite{carrillo2007nonoscillatory,chen2021adaptive,huot2003instability,li2022high,mortezazadeh2024sweep,Qiu_Christlieb,qiu_shu_sl,xiong2019conservative}, finite volume (FV) \cite{ abreu2017new, benkhaldoun2015family,crouseilles2010conservative, healy1993finite, huang2012eulerian, lauritzen2010conservative, nakao2022eulerian, phillips2001conservative, ii2007cip, iske2004conservative}.

Finite volume methods are commonly used to solve hyperbolic conservation laws since they are locally mass conservative due to the numerical flux on the boundaries, able to easily incorporate high spatial resolution methods, and flexible when extending to nonuniform grids. A few recent EL finite volume schemes have used the Rankine-Hugoniot jump condition to define backward-tracing approximate characteristics and partition the space-time region. In particular, Huang and Arbogast proposed an EL finite volume WENO scheme for nonlinear problems in \cite{huang2017eulerian}, and Nakao, Chen, and Qiu developed an EL finite volume method for convection and convection-diffusion equations in \cite{nakao2022eulerian} that used WENO-AO methods \cite{wenoao2016}. Huang, Arbogast, and Hung also proposed a SL finite difference WENO scheme for scalar nonlinear problems in \cite{huang2016semi} that used the Rankine-Hugoniot jump condition to trace the characteristics. In contract to traditional finite volume methods coupled with Runge-Kutta time discretizations, these methods are considerably relaxed with respect to the CFL time-stepping constraint. The numerical tests show that these methods perform well with relatively large time-stepping sizes, but more theoretical justification is still needed. Furthermore, these methods are often less effective when handling nonlinear transport-dominated problems containing shocks. We note, however, that the method in \cite{huang2017eulerian} is able to compute post-shock solutions to Burgers' equation with larger time-stepping sizes. To overcome this challenge, we propose a new high-order EL finite volume method for nonlinear hyperbolic problems that both retains the benefits of existing EL methods and has the ability to capture shocks sharply and accurately while relaxing the CFL condition.

Recently, a cell-merging procedure for a forward-tracing partition of the space-time domain was proposed to handle the shock case when characteristics/partition lines intersect. The first-order scheme is shown to be total-variation-diminishing (TVD) and maximum-principle-preserving (MMP) with a time-stepping size that is twice as large as what the regular CFL condition allows for Burgers' equation; see Theorem 3.1 in \cite{yang2023stability}. In this paper, we extend the first-order scheme in \cite{yang2023stability} to high-order accuracy by coupling with nonuniform ENO \cite{cockburn1998essentially,shu1999high} and WENO-AO \cite{wenoao2016} methods for spatial reconstruction, and high-order strong stability-preserving (SSP) Runge-Kutta methods for the temporal discretization. A distinctive feature of this scheme is the partition of the space-time regions via forward-tracing approximations of the characteristics, and the merging of space-time regions when such forward-tracing characteristics intersect in the case of shocks. This is in contrast to the backward-tracing approximations of the characteristics from \cite{nakao2022eulerian}. Since the allowable time-stepping sizes are larger in the EL/SL framework, approximating forward characteristics explicitly offers advantages in accuracy, efficiency and stability as proved in \cite{yang2023stability}. Dimensional splitting is used to extend the one-dimensional scheme to higher dimensions. Numerical tests show that the proposed method achieves high-order accuracy, and sharply captures shocks and rarefaction waves. As with most EL methods, our primary goal here is to follow the characteristics in order to relax the CFL condition while maintaining high-order accuracy and capturing shocks.

The paper is structured as follows. In Section \ref{sec:ELRKFV_backward}, we review the backward EL-RK-FV method from \cite{nakao2022eulerian}. Section \ref{sec:ELRKFV_forward} introduces the proposed forward EL-RK-FV method by detailing the first-order scheme from \cite{yang2023stability}, and subsequently extending to high-order schemes. Section \ref{sec:numerics} includes numerical tests that demonstrate the performance of the proposed high-order scheme. Conclusions are made in Section \ref{sec:conclusion}. The appendix is included at the end.

\section{Review of the backward EL-RK-FV method}
\label{sec:ELRKFV_backward}

\begin{figure}[h!]
    \centering
    \includegraphics[width=0.4\textwidth]{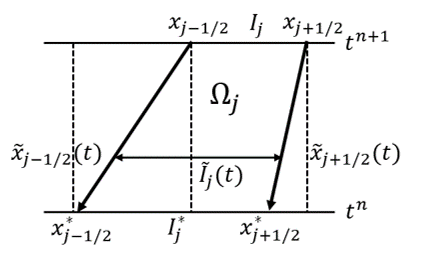}
    \caption{The backward EL-RK-FV space-time region.}
    \label{bkwd_spacetime}
\end{figure}

We start by reviewing a previous EL-RK-FV method that traces the characteristics backward in time \cite{nakao2022eulerian}. This is meant to provide a broad overview of the EL framework of discussion, as well as develop the motivation for our proposed algorithm. Consider the one-dimensional hyperbolic conservation law, 
\begin{equation}
u_t+f(u)_x = 0, \quad x\in(0,2\pi),\quad t>0,
\label{eq: scalar1d}
\end{equation}
where $f(x)$ is a strictly convex flux function, the flow field is a continuous function of space and time, and boundary conditions are periodic. The computational mesh is $0=x_{\frac12}< x_{\frac32}<\cdots< x_{N_x+\frac12} =2\pi$.
Let $I_j=[x_{j-\frac12}, x_{j+\frac12} ]$ denote a cell of length $h_j=x_{j+\frac12}-x_{j-\frac12}$. Further let $t^n$ be the $n$-th time level and $\Delta t =t^{n+1}-t^n$ be the time-stepping size for all $n$. We define the (uniform) cell averages of the solution by
\beq
\bar{u}_j(t) \coloneqq \frac{1}{h_j} \int_{x_{j-\frac12}}^{x_{j+\frac12}} u(x,t)dx.
\label{eq:avg}
\eeq 
Discretizing in time, we solve for $\bar{u}_j^n\approx\bar{u}_j(t^n)$. In the backward EL-RK-FV method and throughout this paper, we use bars $\bar{u}_j$ to denote \textit{uniform} cell averages and tildes $\tilde{u}_j$ to denote \textit{nonuniform} cell averages. First, we partition the space-time region $[0,2\pi]\times [t^n,t^{n+1}]$ into $N_x$ sub-regions by tracing the characteristics backward in time. The velocity of the flow field on the grid back be approximated by the Rankine-Hugoniot jump condition,
\beq
\nu_{j+\frac{1}{2}} = 
\begin{cases}
    \mathlarger{\frac{f(\bar{u}_{j+1}^n) - f(\bar{u}_j^n)}{\bar{u}_{j+1}^n-\bar{u}_j^n}},&\text{if } \bar{u}_j^n\neq\bar{u}_{j+1}^n,\\
    f'(\bar{u}_j^n),&\bar{u}_j^n=\bar{u}_{j+1}^n.
\end{cases}
\label{eq: RHjump}
\eeq
By using the known cell averages at time $t^n$ to approximate the velocity of the characteristics at time $t^{n+1}$, we can trace approximate characteristics backward in time to find the upstream nodes $x_{j+\frac{1}{2}}^*$ for all $j=0,1,...,N_x$. As seen in Figure \ref{bkwd_spacetime}, the space-time region has been partitioned into sub-regions $\Omega_j$, $j=1,2,...,N_x$. Note that the upstream traceback grid is in general nonuniform, and by approximating the characteristics with linear space-time curves we can handle nonlinear hyperbolic problems.

Broadly speaking, the solution is remapped and the uniform cell averages (at time $t^n$) are projected onto the possibly nonuniform traceback grid using WENO-AO reconstruction polynomials \cite{nakao2022eulerian}. The robustness, high-order accuracy, and existence of linear weights at arbitrary points make the WENO-AO method \cite{wenoao2016} attractive for the backward EL-RK-FV scheme. Rewriting equation \eqref{eq: scalar1d} in divergence form $\nabla_{t,x}\cdot(u,f(u))^T=0$, integrating over the space-time region $\Omega_j$, applying the divergence theorem, and writing in the time-differential form, the semi-discrete formulation is
\beq
\frac{d}{dt}\int_{\tilde{I}_j(t)}{u(x,t)dx} = -\left[\hat{F}_{j+\frac{1}{2}}(t)-\hat{F}_{j-\frac{1}{2}}(t)\right],
\label{eq: backward_semidiscrete}
\eeq
where $F_{j+\frac{1}{2}}(t)=f(u(\tilde{x}_{j+\frac{1}{2}}(t),t)) - \nu_{j+\frac{1}{2}}u(\tilde{x}_{j+\frac{1}{2}}(t),t)$ is the modified flux function used to compute an appropriate monotone numerical flux (e.g., Lax-Friedrichs flux),
\[\hat{F}_{j+\frac{1}{2}}(t) = \hat{F}_{j+\frac{1}{2}}(u_{j+\frac{1}{2}}^-,u_{j+\frac{1}{2}}^+;t).\]
Spatial reconstruction of the point values $u(\tilde{x}_{j+\frac{1}{2}}^{\pm}(t),t)$ is done using a WENO-AO method with a nonuniform-to-uniform transformation in which the nonuniform cell averages are used to compute the solution on the linear space-time curves \cite{nakao2022eulerian}. Using a method of lines approach with explicit strong stability-preserving (SSP) Runge-Kutta methods \cite{gottlieb2001strong}, evolve the solution along the approximate characteristics from time $t^n$ to time $t^{n+1}$. Note that the approximate characteristics recover the uniform background grid at time $t^{n+1}$. The backward EL-RK-FV scheme extends to higher dimensions using high-order dimensional splitting methods \cite{forest1990fourth,yoshida1990construction}. We refer the reader to the original paper \cite{nakao2022eulerian} for more details.

The purpose of reviewing the backward EL-RK-FV scheme is to demonstrate the spirit of this class of Eulerian-Lagrangian methods. That is, to use the Rankine-Hugoniot jump condition to partition the space-time region, approximate and project the solution using WENO-AO methods, and evolve along the approximate characteristics using Runge-Kutta methods in a method of lines framework. Doing so maintains high-order accuracy, resolves spurious oscillations due to Gibbs phenomenon, and relaxes the time-stepping constraint from the CFL condition.

Despite the effectiveness of the backward EL-RK-FV method when solving nonlinear hyperbolic conservation laws, there is one major flaw. When capturing shocks (for example, the Riemann problem for Burgers' equation) the time-stepping size is usually constrained by the CFL condition since the linear traceback characteristics will otherwise intersect when finding the upstream nodes. Furthermore, when a shock occurs during a time-step evolution, the partition lines around the shock defined by the Rankine-Hugoniot jump condition (at time $t^n$) are no longer a good approximation to characteristics in the nonlinear shock setting. This defeats the main advantage of the EL framework, namely, the ability to handle large CFL numbers. Hence, we propose a forward-tracing EL-RK-FV scheme that is more suitable for handling intersecting characteristics.

\begin{rem}
Although the backward EL-RK-FV method struggles to handle intersecting characteristics due to shock formations, it can still capture discontinuities so long as the characteristics do not intersect (see the rigid body rotation and swirling deformation examples in \cite{nakao2022eulerian}).
\end{rem}

\section{The forward EL-RK-FV method}
\label{sec:ELRKFV_forward}

Just like the backward EL-RK-FV scheme, a space-time region $\Omega_j$ is defined for the forward EL-RK-FV scheme. We now trace the linear approximate characteristics \textit{forward} in time from $t^n$ to $t^{n+1}$ instead of backward. The linear approximate characteristics at the cell boundaries are still determined by the Rankine-Hugoniot jump condition \eqref{eq: RHjump} with the uniform cell averages $\bar{u}_j^n$.

\begin{figure}[h!]
    \centering
    \includegraphics[width=0.4\textwidth]{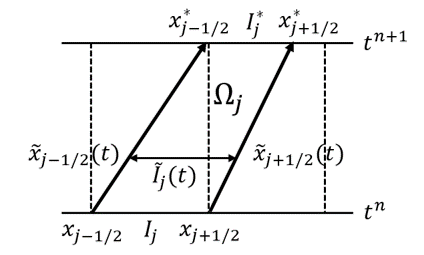}
    \caption{The forward EL-RK-FV space-time region.}
    \label{fwd_spacetime}
\end{figure}

As seen in Figure \ref{fwd_spacetime}, the downstream nodes are defined by $x_{j+\frac{1}{2}}^*\coloneqq x_{j+\frac{1}{2}}+\nu_{j+\frac{1}{2}}\Delta t$, $j=0,1,...,N_x$. The approximate characteristics over the interval $[t^n,t^{n+1}]$ are
\begin{equation}
    \tilde{x}_{j+\frac{1}{2}}(t) = x_{j+\frac{1}{2}} + \nu_{j+\frac{1}{2}}(t-t^n),\qquad t^n\leq t\leq t^{n+1},\qquad j=0,1,...,N_x.
    \label{eq: forward_characteristics}
\end{equation}
Note that $\tilde{x}_{j+\frac{1}{2}}(t^n) = x_{j+\frac{1}{2}}$ and $\tilde{x}_{j+\frac{1}{2}}(t^{n+1}) = x_{j+\frac{1}{2}}^*$. That is, the initial nodes of the characteristics (at time $t^n$) are exactly the uniform background grid; the downstream nodes (at time $t^{n+1}$) form a nonuniform grid in general. The approximate characteristics form the linear space-time boundaries
\[\mathcal{S}_{\ell} = \{(\tilde{x}_{j-\frac{1}{2}}(t),t)\text{ : }t^n\leq t\leq t^{n+1}\} \quad\text{and}\quad \mathcal{S}_{r} = \{(\tilde{x}_{j+\frac{1}{2}}(t),t)\text{ : }t^n\leq t\leq t^{n+1}\}.\]
Let $\tilde{I}_j(t) = [\tilde{x}_{j-\frac{1}{2}}(t),\tilde{x}_{j+\frac{1}{2}}(t)],\ t^n\leq t\leq t^{n+1}$. The space-time region $\Omega_j$ is the region bounded by $I_j$, $I_j^*$, $\mathcal{S}_{\ell}$, and $\mathcal{S}_r$. Assuming the characteristics do not intersect, one could simply evolve the cell averages over $\Omega_j$, $j=1,2,...,N_x$ using a method of lines approach, followed by a projection onto the uniform background mesh. However, partitioning the space-time in the presence of intersecting characteristics due to shock formations requires more care. One can easily confirm if the approximate characteristics intersect by checking if $x_{j+\frac{1}{2}}^*-x_{j-\frac{1}{2}}^*<0$. In which case, the approximate characteristics could intersect at some intermediate time $t$ such that $\tilde{x}_{j+\frac{1}{2}}(t)=\tilde{x}_{j-\frac{1}{2}}(t)$.

We introduce a merging procedure for partitioning the space-time region when dealing with shock formations. This merging procedure for the forward EL-RK-FV scheme was first proposed in \cite{yang2023stability}, where the first-order scheme was proved to be total-variation-diminishing and maximum-principle-preserving for Burgers' equation. Although the stability analysis in \cite{yang2023stability} is only shown for the first-order scheme applied to Burgers' equation, the first-order algorithm and theoretical results provide a solid theoretical foundation for higher-order extensions and general strictly convex flux functions.

The focus of this paper is to extend the first-order algorithm in \cite{yang2023stability} to higher-order schemes. The theoretical proof for the high-order methods is quite complicated. However, our numerical tests suggest that the algorithm performs with good accuracy and stability. We only consider Burgers' equation as an initial step, but the proposed scheme could be generalized to other scalar nonlinear flux functions.

The forward EL-RK-FV scheme presented in this paper has the following major steps:
\begin{enumerate}
	\item Partition the space-time region into sub-regions using a merging procedure so that the (approximate) characteristics propagating from the cell boundaries at time $t^n$ do not intersect within a single time-step.
	\item Evolve the cell averages over the newly merged cells from time $t^n$ to $t^{n+1}$ using a high-order strong stability-preserving Runge-Kutta method \cite{gottlieb2001strong}. Nonuniform ENO and WENO-AO schemes are used for spatial reconstruction.
	\item Project the solution at time $t^{n+1}$ over the nonuniform grid onto the uniform background grid using nonuniform reconstruction.
\end{enumerate}


\subsection{A first-order scheme for Burgers' equation}

We start by outlining the first-order scheme previously shown to be total-variation-diminishing and maximum-principle-preserving for Burgers' equation \cite{yang2023stability}. For the following definitions, we let $[0,2\pi]=\bigcup\limits_{j=1}^{N_x}{I_j}$, where $I_j$ are ordered disjoint uniform cells.

\begin{defn}\label{defn: troubledcell}
A cell $I_j$ is said to be troubled if either: (1) the approximate characteristics \eqref{eq: forward_characteristics} originating from the cell boundaries $x_{j\pm\frac{1}{2}}$ intersect at some time $t\in(t^n,t^{n+1})$, or (2) $\bar{u}_j^n$ is significantly larger than $\bar{u}_{j+1}^n$ or significantly smaller than $\bar{u}_{j-1}^n$. Further let $\lambda>0$ be an upper bound for the time-stepping constraint such that $\Delta t\leq\lambda\Delta x$. A troubled cell falls into one of five mutually exclusive classifications:
\begin{itemize}
	\item[1.] $I_j$ is called a troubled cell of type I for $\lambda$ if $\mathlarger{f'(\bar{u}_{j-1}^n)>f'(\bar{u}_{j+1}^n) + \frac{2}{\lambda}}$. Since $f'(u)=u$ for Burgers' equation, we check if $\mathlarger{\bar{u}_{j-1}^n>\bar{u}_{j+1}^n+\frac{2}{\lambda}}$.
	\item[2.] $I_j$ is called a troubled cell of type II for $\lambda$ if it is not a troubled cell of type I, and satisfies $\mathlarger{\bar{u}_{j-1}^n>\bar{u}_j^n+\frac{2}{\lambda}}$ and $\bar{u}_{j-1}^n\geq\bar{u}_{j+1}^n\geq\bar{u}_j^n$.
	\item[3.] $I_j$ is called a troubled cell of type III for $\lambda$ if it is not a troubled cell of type I, and satisfies $\mathlarger{\bar{u}_j^n>\bar{u}_{j+1}^n+\frac{2}{\lambda}}$ and $\bar{u}_j^n\geq\bar{u}_{j-1}^n\geq\bar{u}_{j+1}^n$.
	\item[4.] $I_j$ is called a troubled cell of type IV for $\lambda$ if it is neither a troubled cell of type I nor type II, and satisfies $\mathlarger{\bar{u}_{j-1}^n>\bar{u}_j^n+\frac{2}{\lambda}}$.
	\item[5.] $I_j$ is called a troubled cell of type V for $\lambda$ if it is neither a troubled cell of type I nor type III, and satisfies $\mathlarger{\bar{u}_j^n>\bar{u}_{j+1}^n+\frac{2}{\lambda}}$.
\end{itemize}
\end{defn}

The conditions for the five types of troubled cells in Definition \ref{defn: troubledcell} are specific to Burgers' equation since $f'(u)=u$. A troubled cell of type I occurs if the approximate characteristics intersect at some time $t\in(t^n,t^{n+1})\subseteq(t^n,t^n+\lambda\Delta x)$ since
\[0>x_{j+\frac{1}{2}}^*-x_{j-\frac{1}{2}}^*=\Delta x+\Delta t(\nu_{j+\frac{1}{2}}-\nu_{j-\frac{1}{2}})=\Delta x+\frac{\Delta t}{2}(\bar{u}_{j+1}^n-\bar{u}_{j-1}^n),\]
and hence
\begin{equation}\label{eq: troubledcell_type1}
\bar{u}_{j-1}^n>\bar{u}_{j+1}^n+\frac{2\Delta x}{\Delta t}\geq\bar{u}_{j+1}^n+\frac{2}{\lambda}.
\end{equation}

\begin{figure}[h!]
	\centering
	\includegraphics[width=0.9\textwidth]{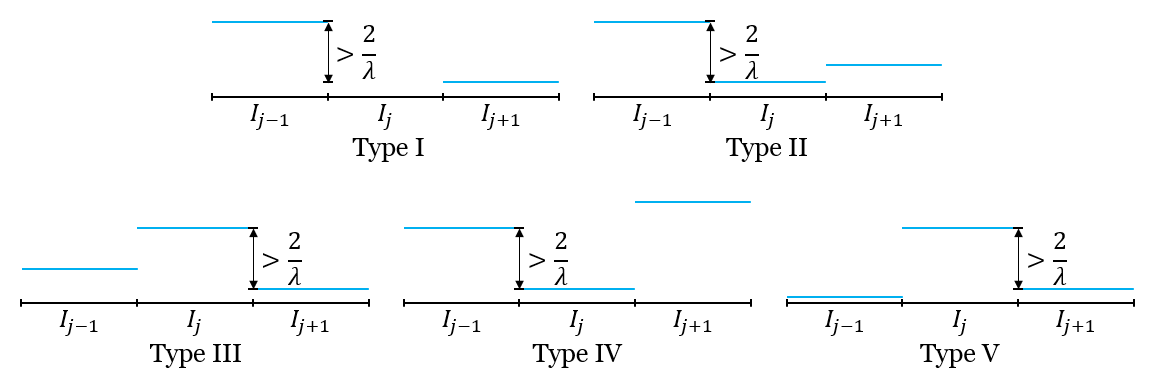}
	\caption{A visualization of the five types of troubled cells defined in Definition \ref{defn: troubledcell}.}
	\label{fig: troubledcells}
\end{figure}

Observe that equation \eqref{eq: troubledcell_type1} suggests a strong shock occurs when $\bar{u}_{j-1}^n-\bar{u}_{j+1}^n>2/\lambda$. Troubled cells of types II-V account for the other instances in which similar strong shocks may occur. In particular, troubled cells of types II and III occur in the presence of strong shock waves propagating from the cell boundaries. Whereas, troubled cells of types IV and V are unlikely but can occur in the presence of rarefaction waves with strong numerical oscillations. Figure \ref{fig: troubledcells} provides a visualization of the five types of troubled cells.

\begin{figure}[h!]
    \centering
    \includegraphics[width=0.9\textwidth]{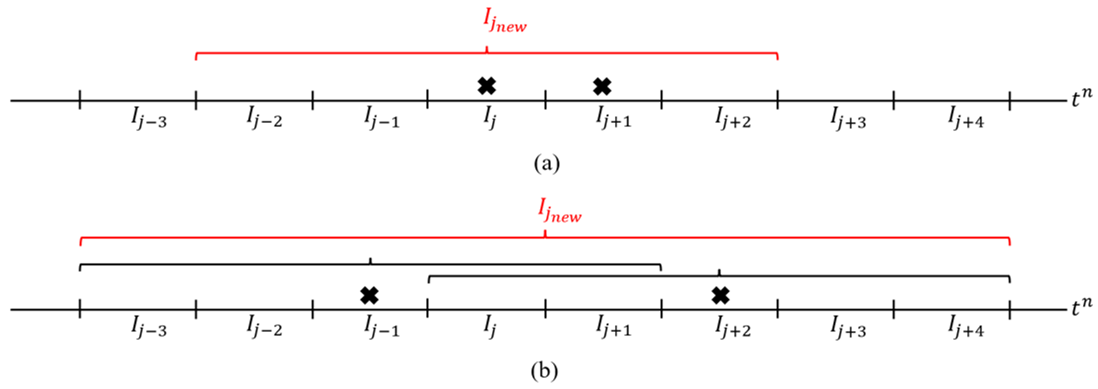}
    \caption{(a) Two troubled cells over the stencil $\{I_{j-1},I_j,I_{j+1}\}$ for which the effective troubled cell that determines the influence region is $I_j$. (b) Step 3, the merging of two overlapping influence regions from two neighboring effective troubled cells are not isolated.}
    \label{fig: ETC}
\end{figure}

Since the characteristics could quickly intersect near shocks and hence severely limit the allowable time-stepping size, we propose merging a collection of cells surrounding each troubled cell. In doing so, we bypass the intersecting characteristics and take larger time-steps. Moving from left to right: if there are multiple troubled cells over a stencil $\{j-1,j,j+1\}$, then we merge the cells only based on the left-most cell; we define this as the \textit{effective troubled cell}. One could also define the effective troubled cell based on the right-most cell moving right to left; it does not matter as long as one is consistent throughout. Figure \ref{fig: ETC}(a) shows an effective trouble cell based on the $I_j$ instead of $I_{j+1}$. If there is only a single troubled cell over the stencil $\{j-1,j,j+1\}$, then we consider it an isolated effective troubled cell.

Our merging procedure goes as follows. Identify the effective troubled cells; and for each effective troubled cell, merge a local collection of cells together based on the type of troubled cell. We define the so-called \textit{influence region} of each effective troubled cell that will determine how we merge the cells.

\begin{defn}\label{defn: influenceregion}
For each effective troubled cell $I_j$, we consider the cells $I_i$, $i=j-3,...,j+3$. Assume the initial condition is bounded by $b\leq u_0(x)\leq a$. The influence region of the effective troubled cell $I_j$ is defined as follows:
\begin{itemize}
	\item[1.] If $I_j$ is a troubled cell of type IV, then $I_{j-1}$ is a troubled cell of type V. Moreover, the influence region contains cells $I_i$, $i=j-2,...,j+1$.
	\item[2.] Assume $\mathlarger{\bar{u}_{j-1}^n+\bar{u}_{j}^n+\bar{u}_{j+1}^n>\frac{7a+5b}{4}}$. If $\mathlarger{\bar{u}_{j+2}^n<\frac{a+3b}{4}}$ or $\mathlarger{\bar{u}_{j+2}^n+\bar{u}_{j+3}^n<\frac{a+3b}{2}}$, then the influence region contains $I_i$, $i=j-2,...,j+3$.
	\item[3.] Assume $\mathlarger{\bar{u}_{j-1}^n+\bar{u}_{j}^n+\bar{u}_{j+1}^n>\frac{5a+7b}{4}}$. If $\mathlarger{\bar{u}_{j-2}^n>\frac{3a+b}{4}}$ or $\mathlarger{\bar{u}_{j-3}^n+\bar{u}_{j-2}^n>\frac{3a+b}{2}}$, then the influence region contains $I_i$, $i=j-3,...,j+2$.
	\item[4.] In all other cases, the influence region contains $I_i$, $i=j-2,...,j+2$.
\end{itemize}
Two distinct effective troubled cells are said to be \textit{isolated} if their influence regions do not overlap.
\end{defn}

The influence region for most effective troubled cells will contain cells $I_i$, $i=j-2,...,j+2$, that is, case 4. Cases 1-3 in Definition \ref{defn: influenceregion} account for particular situations when the influence region might need one more cell or one less cell on either the left or right. Case 1 in Definition \ref{defn: influenceregion} applies to rarefaction waves. Whereas, cases 2 and 3 in Definition \ref{defn: influenceregion} apply in the special and rare situation of a shock-shock interaction where the intersection of characteristics occurs very close to the next time level $t^{n+1}$. In such cases, we must include more information from either the left or right of the neighboring trouble cells. We note that the specific conditions in Definition \ref{defn: influenceregion} are specially determined so that the first-order scheme is total-variation-diminishing and maximum-principle-preserving \cite{yang2023stability}.

\begin{figure}[h!]
    \centering
    \includegraphics[width=0.6\textwidth]{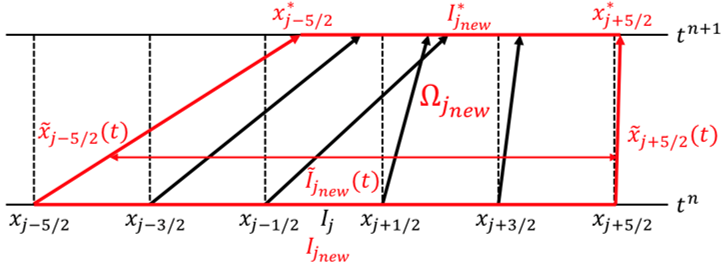}
    \caption{The merging of the cells in the influence region of an isolated effective troubled cell; see case 4 in Definition \ref{defn: influenceregion}.}
    \label{fig: mergingexample}
\end{figure}

For each effective troubled cell, we merge the cells in their influence region. If two neighboring effective troubled cells are not isolated (i.e., their influence regions overlap), then we merge the cells from both influence regions; Figure \ref{fig: ETC}(b) shows such an example. We denote the cells that make up the post-merging grid by $I_{j_{\text{new}}}$, $j_{\text{new}}=1,2,...,N_{\text{merged}}$. Figure \ref{fig: mergingexample} illustrates the merging of five cells in the influence region of an isolated effective troubled cell. Note that the approximate characteristics of the post-merging grid do not intersect.

For each influence region, the cell averages over the influence region are averaged. Figure \ref{fig: mergingexample} provides an illustrative example for which the cell average over $I_{j_{\text{new}}}$ at time $t^n$ would be
\begin{equation}
	\tilde{u}_{j_{\text{new}}}^n = \frac{1}{\Delta x_{j_{\text{new}}}}\int_{I_{j_{\text{new}}}}{u(x,t^n)dx} \approx \frac{\bar{u}_{j-2}^n + \bar{u}_{j-1}^n + \bar{u}_j^n + \bar{u}_{j+1}^n + \bar{u}_{j+2}^n}{5},
\end{equation}
or rather,
\begin{equation}
	\int_{I_{j_{\text{new}}}}{u(x,t^n)dx} \approx \Delta x\left(\bar{u}_{j-2}^n + \bar{u}_{j-1}^n + \bar{u}_j^n + \bar{u}_{j+1}^n + \bar{u}_{j+2}^n\right).
\end{equation}

\begin{rem}\label{rem: timesteprestriction}
The stability proof in \cite{yang2023stability} showing the first-order forward EL-RK-FV scheme is total-variation-diminishing and maximum-principle-preserving assumes the time-stepping constraint
\begin{equation}\label{eq: timeconstraint}
	\Delta t<\lambda\Delta x,\qquad\text{with}\qquad \lambda=\frac{4}{a-b},
\end{equation}
where $b\leq u_0(x)\leq a$. The allowable time-stepping sizes from constraint \eqref{eq: timeconstraint} are at least twice as large as the standard Courant-Friedrichs-Lewy (CFL) condition. Alternatively, the stability proof also holds if $a$ and $b$ are chosen to be the local maximum and minimum in the influence region of the effective troubled cell plus the two neghboring cells of the influence region \cite{yang2023stability}. For example, if the influence region of the effective troubled cell contains $I_i$, $i=j-2,...,j+2$, then
\begin{equation}
	\Delta t<\lambda\Delta x,\qquad\text{with}\qquad \lambda=\frac{4}{\max\limits_{j}{\{a_j-b_j\}}},
\end{equation}
where
\[a_j=\max\limits_{i=0,\pm 1,\pm 2,\pm 3}{\bar{u}_{j+i}^n},\qquad b_j=\min\limits_{i=0,\pm 1,\pm 2,\pm 3}{\bar{u}_{j+i}^n}.\]
In all of our numerical experiments, we use the global time-stepping constraint in equation \eqref{eq: timeconstraint}.
\end{rem}

\begin{rem}
In general, the influence region of a troubled cell $I_j$ includes cells $I_i$, $i=j-m,...,j+m$ (plus or minus one cell from the left or right in specific cases), for some integer $m\geq 2$. Definition \ref{defn: influenceregion} assumes $m=2$ since it corresponds to the smallest influence region that guarantees stability. Larger $m$ can be taken, but the solution will be averaged over more cells and could negatively affect the resolution of the solution near effective troubled cells that involve shocks.
\end{rem}

Now that we have introduced a merging procedure in Definitions \ref{defn: troubledcell} and \ref{defn: influenceregion} for handling intersecting characteristics and strong physical interactions, we can set up the scheme formulation for evolving the solution. As seen in Figure \ref{fig: mergingexample}, the post-merging grid $\{I_{j_{\text{new}}}\text{ : }j_{\text{new}}=1,2,...,N_{\text{merged}}\}$ and the approximate characteristics $\tilde{x}_{j_{\text{new}}\pm\frac{1}{2}}(t)$ determined by equation \eqref{eq: forward_characteristics} define the space-time regions $\Omega_{j_{\text{new}}}$, $j_{\text{new}}=1,2,...,N_{\text{merged}}$. Similar to the backward EL-RK-FV method, we integrate over the space-time region $\Omega_{j_{\text{new}}}$ and rewrite the system into a time-differential form to get
\beq
\frac{d}{dt}\int_{\tilde{I}_{j_{\text{new}}}(t)}{u(x,t)dx} = -\left[\hat{F}_{j_{\text{new}}+\frac{1}{2}}(t)-\hat{F}_{j_{\text{new}}-\frac{1}{2}}(t)\right],
\label{eq: forward_semidiscrete}
\eeq
where the modified flux function is $F_{j_{\text{new}}+\frac{1}{2}}(t)=f(u(\tilde{x}_{j_{\text{new}}+\frac{1}{2}}(t),t)) - \nu_{j_{\text{new}}+\frac{1}{2}}u(\tilde{x}_{j_{\text{new}}+\frac{1}{2}}(t),t)$, and $\hat{F}_{j_{\text{new}}+\frac{1}{2}}(t) = \hat{F}_{j_{\text{new}}+\frac{1}{2}}(u_{j_{\text{new}}+\frac{1}{2}}^-,u_{j_{\text{new}}+\frac{1}{2}}^+;t)$ still denotes any appropriate monotone numerical flux. We use Lax-Friedrichs flux throughout this paper. All that remains is to update equation \eqref{eq: forward_semidiscrete} from $t^n$ to $t^{n+1}$ as a method of lines system. Evolving equation \eqref{eq: forward_semidiscrete} with the forward Euler method,
\begin{equation}
    \int_{I_{j_{\text{new}}}^*}{u(x,t^{n+1})dx} = \int_{I_{j_{\text{new}}}}{u(x,t^n)dx} - \Delta t\left[\hat{F}_{j_{\text{new}}+\frac{1}{2}}(t^n) - \hat{F}_{j_{\text{new}}-\frac{1}{2}}(t^n)\right].
    \label{eq: forwardEuler}
\end{equation}

The updated uniform cell averages $\bar{u}_j^{n+1}$ are obtained by projecting the updated non-uniform cell averages $\tilde{u}_{j_{\text{new}}}^{n+1} = \int_{I_{j_{\text{new}}}^*}{u(x,t^{n+1})dx}/\Delta x_{j_{\text{new}}}^*$ onto the background uniform mesh. Since this is just a first-order scheme, the $L^2$ projection onto the backgroud uniform mesh is
\begin{equation}\label{eq: firstorder_projection}
\int\limits_{I_j}{u(x,t^{n+1})dx} \approx
\begin{cases}
	\mathlarger{\frac{\Delta x}{\Delta x_{j_{\text{new}}}^*}\int\limits_{I_{j_{\text{new}}}^*}{u(x,t^{n+1})dx}},\hspace*{5cm}\text{if}\ I_j\subseteq I_{j_{\text{new}}}^*,\\
	\ \\
	\mathlarger{\frac{x_{\ell+\frac{1}{2}}^* - x_{j-\frac{1}{2}}}{\Delta x_{\ell}^*}\int\limits_{I_{{\ell}}^*}{u(x,t^{n+1})dx} + \sum\limits_{j_{\text{new}}=\ell+1}^{r-1}{\int\limits_{I_{j_{\text{new}}}^*}{u(x,t^{n+1})dx}} + \frac{x_{j+\frac{1}{2}} - x_{r-\frac{1}{2}}^*}{\Delta x_{r}^*}\int\limits_{I_{{r}}^*}{u(x,t^{n+1})dx}},\\
	\hspace*{9cm}\text{if}\ I_j\subseteq \bigcup\limits_{j_{\text{new}}=\ell}^{r}{I_{j_{\text{new}}}^*}\ \text{for $\ell<r$}.
\end{cases}
\end{equation}

Prior to merging, the uniform cell averages $\bar{u}_j^n$ are used in the first stage of the RK method for first-order spatial reconstruction in the numerical flux at time $t^n$. The first-order forward EL-RK-FV scheme is given by Algorithm \ref{algo1}.

\begin{algorithm}[t!]
	\caption{First-order forward EL-RK-FV method}
	\label{algo1}
		{\bf Input:}  $\bar{u}^{n}_j$, $j=1,2,\dots, N_x $ \\
	    {\bf Output:} $\bar{u}^{n+1}_j$, $j=1,2,\dots, N_x$  
	\begin{algorithmic}[1]
	\For{$j=0,1,\dots, N_x$}
		\State Compute velocities $v_{j+\frac12}$ using equation \eqref{eq: RHjump}.
		\State Locate the downstream nodes $x_{j+\frac{1}{2}}^*$ using equation \eqref{eq: forward_characteristics}.
	\EndFor
		
	\For{$j=0,1,\dots, N_x $}
		 \State Compute $\hat{F}_{j+ \frac12}(t^n)$ at the (uniform) cell boundary using $u_{j+\frac{1}{2}}^{n,-}\approx\overline{u}_j^n$ and $u_{j+\frac{1}{2}}^{n,+}\approx\overline{u}_{j+1}^n$.
	\EndFor

	\For{$j=1,2,\dots, N_x$}
	    \State Determine if $I_j$ is a troubled cell of type I-V using Definition \ref{defn: troubledcell}.
	    \State According to Definition \ref{defn: troubledcell}, keep only the effective troubled cells based on the left-most troubled cell.
	\EndFor
	    
	 \For{each effective troubled cell}
	    \State Determine the influence region according to Definition \ref{defn: influenceregion}, and merge the cells accordingly.
	    \State If neighboring influence regions overlap, then merge the cells from both influence regions.
	 \EndFor
	 \State We now have the post-merging grid (at time $t^n$) $\{I_{j_{\text{new}}}\text{ : }j_{\text{new}}=1,2,...,N_{\text{merged}}\}$.
	 
	 \For{$j_{\text{new}}=1,2,\dots,N_{\text{merged}}$}
	    \State $I_{j_{\text{new}}}$ consists of the uniform cells $I_i$ for $i\in\mathcal{S}_{j_{\text{\text{new}}}}$.
		\State Gather the information for $I_{j_{new}}$:
		\State\qquad Take $\nu_{j_{\text{new}}\pm\frac{1}{2}}$ from $\{\nu_{j+\frac{1}{2}}\text{ : }j=0,1,...,N_x\}$.
		\State\qquad Take $\hat{F}_{j_{\text{new}}\pm\frac{1}{2}}(t^n)$ from $\{\hat{F}_{j+\frac{1}{2}}(t^n)\text{ : }j=0,1,...,N_x\}$.
		\State\qquad Compute $\int_{I_{j_{\text{new}}}}{u(x,t^n)dx}\approx\Delta x\sum\limits_{i\in\mathcal{S}_{j_{\text{new}}}}{\bar{u}_i^n}$.
	    \State Compute $\int_{I_{j_{\text{new}}}^*}{u(x,t^{n+1})dx}$ using equation \eqref{eq: forwardEuler}, and hence obtain $\tilde{u}_{j_{\text{new}}}^{n+1}$.
	 \EndFor

	 \For{$j=1,2,\dots, N_x$}
		\State We have that $I_j\subseteq\bigcup\limits_{j_{\text{new}}=\ell}^{j_{\text{new}}=r}{I_{j_{\text{new}}}^*}$ for some $\ell\leq r$.
		\If{$\ell=r$}
		    \State Compute $\bar{u}_j^{n+1}\approx\tilde{u}_{\ell}^{n+1}$.
		\Else
		    \State Compute $\mathlarger{\bar{u}_j^{n+1}\approx \frac{1}{\Delta x}\left((x_{\ell+\frac{1}{2}}^*-x_{j-\frac{1}{2}})\tilde{u}_{\ell}^{n+1} + \sum\limits_{j_{\text{new}}=\ell+1}^{j_{\text{new}}=r-1}{\Delta x_{j_{\text{\text{new}}}}^*\tilde{u}_{j_{\text{new}}}^{n+1}} +  (x_{j+\frac{1}{2}} - x_{r-\frac{1}{2}}^*)\tilde{u}_{r}^{n+1}\right)}$.
		 \EndIf
	 \EndFor
	\end{algorithmic}
\end{algorithm}


\subsection{Higher-order schemes for Burgers' equation}

There is a key condition that one must account for when using high-order Runge-Kutta methods in the forward EL-RK-FV scheme. The effective troubled cells defined in Definition \ref{defn: troubledcell} should be identified by checking over the time interval $(t^n,t^{n+2})$, that is, over two time-steps $2\Delta t$. The inequalities in Definition \ref{defn: troubledcell} correspond to checking over a single time-step; see equation \eqref{eq: troubledcell_type1}. Checking over two time-steps when defining the effective troubled cells is justified by looking at the time discretizations. Written in a more convenient form, the optimal second- and third-order SSP RK methods \cite{gottlieb2001strong} for solving the initial value problem $\mathlarger{\frac{dU}{dt}=\mathcal{L}(U)}$, $U(t^0)=U^0$ are respectively
\begin{subequations}
\begin{align}
\begin{split}
    U^{(1)} &= U^n + \Delta t\mathcal{L}(U^n),
\end{split}\\
\begin{split}\label{eq: SSPRK2_t2}
    U^{n+1} &= \frac{1}{2}U^n + \frac{1}{2}\Big(U^{(1)} + \Delta t\mathcal{L}(U^{(1)})\Big),
\end{split}
\end{align}
\label{eq: SSPRK2_v1}
\end{subequations}
and
\begin{subequations}
\begin{align}
\begin{split}
    U^{(1)} &= U^n + \Delta t\mathcal{L}(U^n),
\end{split}\\
\begin{split}\label{eq: SSPRK3_t2}
    U^{(2)} &= \frac{3}{4}U^n + \frac{1}{4}\Big(U^{(1)} + \Delta t\mathcal{L}(U^{(1)})\Big),
\end{split}\\
\begin{split}\label{eq: SSPRK3_t3}
    U^{n+1} &= \frac{1}{3}U^n + \frac{2}{3}\Big(U^{(2)} + \Delta t\mathcal{L}(U^{(2)})\Big).
\end{split}
\end{align}
\label{eq: SSPRK3_v1}
\end{subequations}

Looking at equation \eqref{eq: SSPRK2_t2}, the second stage of the optimal second-order SSP RK method is at time $t = t^{n+1} = \frac{1}{2}(t^n) + \frac{1}{2}(t^n+2\Delta t)$. That is, the updated solution takes average of $U^n$ and a forward Euler approximation using $U^{(1)}$. Similarly, the second stage of the optimal third-order SSP RK method in equation \eqref{eq: SSPRK3_t2} is at time $t = t^{(2)} = t^{n+1/2} = \frac{3}{4}(t^n) + \frac{1}{4}(t^n + 2\Delta t)$; it takes a weighted average of $U^n$ and a forward Euler approximation using $U^{(1)}$. The third stage of the optimal third-order SSP RK method in equation \eqref{eq: SSPRK3_t3} is at time $t = t^{n+1} = \frac{1}{3}(t^n) + \frac{2}{3}(t^{n} + \frac{3}{2}\Delta t)$; it takes a weighted average of $U^n$ and a forward Euler approximation using $U^{(2)}$. Since both time discretizations \eqref{eq: SSPRK2_v1} and \eqref{eq: SSPRK3_v1} rely on approximations of the solutions up to time $t^{n+2}$, we must look over the time interval $(t^n,t^{n}+2\Delta t)\subseteq(t^n,t^n+2\lambda \Delta x)$ when defining the effective troubled cells. Letting $x_{j\pm\frac{1}{2}}^{**}\coloneqq x_{j\pm\frac{1}{2}} + 2\nu_{j\pm\frac{1}{2}}\Delta t$, we must check if
\[0 > x_{j+\frac{1}{2}}^{**} - x_{j-\frac{1}{2}}^{**} = \Delta x + 2\Delta t(\nu_{j+\frac{1}{2}} - \nu_{j-\frac{1}{2}}) = \Delta x + \Delta t(\bar{u}_{j+1}^n - \bar{u}_{j-1}),\]
and hence
\begin{equation}
	\bar{u}_{j-1}^n > \bar{u}_{j+1}^n + \frac{\Delta x}{\Delta t} \geq \bar{u}_{j+1}^n + \frac{1}{\lambda}.
\end{equation}

Identifying the other four types of troubled cells over $(t^n,t^{n+2})$ follows similarly. Although the (effective) troubled cells are determined by looking over $(t^n,t^{n+2})$, the merging procedure and influence regions are still performed over $(t^n,t^{n+1})$ when defining the space-time regions $\Omega_{j_{\text{new}}}$.

\begin{defn}\label{defn: troubledcell_higherorder}
For the second- and third-order forward EL-RK-FV schemes, a troubled cell falls into one of the same five mutually exclusive classifications in Definition \ref{defn: troubledcell}, but replace $\frac{2}{\lambda}$ in the conditions with $\frac{1}{\lambda}$.
\end{defn}

Solving equation \eqref{eq: forward_semidiscrete} using the optimal second-order SSP RK method in equation \eqref{eq: SSPRK2_v1} but expressed as in a Butcher table,

\begin{subequations}
\begin{align}
\begin{split}\label{eq: SSPRK2_t1_v2}
    \int_{\tilde{I}_{j_{\text{new}}}(t^{(1)})}{u(x,t^{(1)})dx} &= \int_{I_{j_{\text{new}}}}{u(x,t^n)dx} - \Delta t\left[\hat{F}_{j_{\text{new}}+\frac{1}{2}}(t^n) - \hat{F}_{j_{\text{new}}-\frac{1}{2}}(t^n)\right].
\end{split}\\
\begin{split}\label{eq: SSPRK2_t2_v2}
    \int_{I_{j_{\text{new}}}^*}{u(x,t^{n+1})dx} &= \int_{I_{j_{\text{new}}}}{u(x,t^n)dx} - \frac{\Delta t}{2}\left[\hat{F}_{j_{\text{new}}+\frac{1}{2}}(t^n) - \hat{F}_{j_{\text{new}}-\frac{1}{2}}(t^n)\right]\\
    &\qquad\qquad\qquad\qquad\qquad - \frac{\Delta t}{2}\left[\hat{F}_{j_{\text{new}}+\frac{1}{2}}(t^{(1)}) - \hat{F}_{j_{\text{new}}-\frac{1}{2}}(t^{(1)})\right],
\end{split}
\end{align}
\label{eq: SSPRK2_v2}
\end{subequations}
\ \\
Similarly, the optimal third-order SSP RK method in equation \eqref{eq: SSPRK3_v1} can be used to solve equation \eqref{eq: forward_semidiscrete}.

\begin{subequations}
\begin{align}
\begin{split}
    \int_{\tilde{I}_{j_{\text{new}}}(t^{(1)})}{u(x,t^{(1)})dx} &= \int_{I_{j_{\text{new}}}}{u(x,t^n)dx} - \Delta t\left[\hat{F}_{j_{\text{new}}+\frac{1}{2}}(t^n) - \hat{F}_{j_{\text{new}}-\frac{1}{2}}(t^n)\right].
\end{split}\\
\begin{split}
    \int_{\tilde{I}_{j_{\text{new}}}(t^{(2)})}{u(x,t^{(2)})dx} &= \int_{I_{j_{\text{new}}}}{u(x,t^n)dx} - \frac{\Delta t}{4} \left[\hat{F}_{j_{\text{new}}+\frac{1}{2}}(t^n) - \hat{F}_{j_{\text{new}}-\frac{1}{2}}(t^n)\right]\\
    &\qquad\qquad\qquad\qquad\qquad - \frac{\Delta t}{4}\left[\hat{F}_{j_{\text{new}}+\frac{1}{2}}(t^{(1)}) - \hat{F}_{j_{\text{new}}-\frac{1}{2}}(t^{(1)})\right],
\end{split}\\
\begin{split}
    \int_{I_{j_{\text{new}}}^*}{u(x,t^{n+1})dx} &= \int_{I_{j_{\text{new}}}}{u(x,t^n)dx} - \frac{\Delta t}{6}\left[\hat{F}_{j_{\text{new}}+\frac{1}{2}}(t^n) - \hat{F}_{j_{\text{new}}-\frac{1}{2}}(t^n)\right]\\
    &\qquad\qquad\qquad\qquad\qquad - \frac{2\Delta t}{3}\left[\hat{F}_{j_{\text{new}}+\frac{1}{2}}(t^{(2)}) - \hat{F}_{j_{\text{new}}-\frac{1}{2}}(t^{(2)})\right]\\
    &\qquad\qquad\qquad\qquad\qquad\qquad - \frac{\Delta t}{6}\left[\hat{F}_{j_{\text{new}}+\frac{1}{2}}(t^{(1)}) - \hat{F}_{j_{\text{new}}-\frac{1}{2}}(t^{(1)})\right],
\end{split}
\end{align}
\label{eq: SSPRK3_v2}
\end{subequations}

High-order ENO and WENO-AO methods are used for spatial reconstruction of the point values $u(\tilde{x}_{j_{\text{new}}+\frac12}^{+}(t),t)$, as well as projecting the nonuniform cell averages (at time $t^{n+1}$) onto the background uniform mesh. In our numerical tests, we use third-order spatial reconstruction for both the second- and third-order algorithms; higher-order reconstructions can also be used. To help maintain brevity, the nonuniform ENO and WENO-AO reconstructions are included in the appendix.

We note that \textit{uniform} WENO-AO \cite{wenoao2016} is used for spatial reconstruction at time $t^n$ since the pre-merging upstream nodes align with the background uniform grid, as seen in Figure \ref{fwd_spacetime}. These values are stored and used for the point approximations over the post-merging grid at time $t^n$, similar to what was done in Algorithm \ref{algo1}.

Similar to the projection in equation \eqref{eq: firstorder_projection} for the first-order scheme, we use high-order ENO or WENO-AO reconstruction polynomials to project the nonuniform cell averages (at time $t^{n+1}$) onto the background uniform mesh to obtain the desired solution $\bar{u}_j^{n+1}$. Let $P_{j_{\text{new}}}(x;t^{n+1})$ be the high-order nonuniform ENO or WENO-AO reconstruction polynomial for cell $I_{j_{\text{new}}}$ at time $t^{n+1}$. Since the linear weights in these reconstruction polynomials exist at arbitrary points, we can integrate over subsets of each cell. As such, the $L^2$ projection onto the background uniform mesh is

\begin{equation}\label{eq: highorder_projection}
\int\limits_{I_j}{u(x,t^{n+1})dx} \approx
\begin{cases}
	\mathlarger{\int\limits_{I_j}{P_{j_{\text{new}}}(x;t^{n+1})dx}},\hspace*{4.5cm}\text{if}\ I_j\subseteq I_{j_{\text{new}}}^*,\\
	\ \\
	\mathlarger{\int\limits_{I_j\cap I^*_{\ell}}{P_{\ell}(x;t^{n+1})dx} + \sum\limits_{j_{\text{new}}=\ell+1}^{r-1}{\int\limits_{I^*_{j_{\text{new}}}}{P_{j_{\text{new}}}(x;t^{n+1})dx}} + \int\limits_{I_j\cap I^*_{r}}{P_{r}(x;t^{n+1})dx}},\\
	\hspace*{7.5cm}\text{if}\ I_j\subseteq \bigcup\limits_{j_{\text{new}}=\ell}^{r}{I_{j_{\text{new}}}^*}\ \text{for $\ell<r$}.
\end{cases}
\end{equation}

Note that equation \eqref{eq: highorder_projection} reduces to equation \eqref{eq: firstorder_projection} if $P_{j_{\text{new}}}(x;t^{n+1})=\tilde{u}_{j_{\text{new}}}^{n+1}$. We may also use the nonuniform reconstruction polynomials $P_{j_{\text{new}}}(x;t)$ using the nonuniform cell averages at time $t$ to compute the numerical flux $\hat{F}_{j_{\text{new}}\pm\frac12}(t)$. The second-order forward EL-RK-FV scheme (with third-order spatial reconstruction) is given by Algorithm \ref{algo2}. Higher-order schemes follow similarly.

\begin{algorithm}[t!]
	\caption{Second-order forward EL-RK-FV method}
	\label{algo2}
		{\bf Input:}  $\bar{u}^{n}_j$, $j=1,2,\dots, N_x $ \\
	    {\bf Output:} $\bar{u}^{n+1}_j$, $j=1,2,\dots, N_x$  
	\begin{algorithmic}[1]
	\For{$j=0,1,\dots, N_x$}
		\State Compute velocities $v_{j+\frac12}$ using equation \eqref{eq: RHjump}.
		\State Locate the downstream nodes $x_{j+\frac{1}{2}}^*$ and $\tilde{x}_{j+\frac{1}{2}}^{**}$.
	\EndFor
	
	\For{$j=0,1,\dots, N_x $}
		 \State Compute $\hat{F}_{j+ \frac12}(t^n)$ at the (uniform) cell boundary using \textit{uniform} WENO-AO3 \cite{wenoao2016} to compute $u_{j+\frac{1}{2}}^{n,\pm}$.
	\EndFor
	
	\For{$j=1,2,\dots, N_x$}
	    \State Determine if $I_j$ is a troubled cell of type I-V using Definition \ref{defn: troubledcell_higherorder}.
	    \State According to Definition \ref{defn: troubledcell_higherorder}, keep only the effective troubled cells based on the left-most troubled cell.
	\EndFor
	    
	 \For{each effective troubled cell}
	    \State Determine the influence region according to Definition \ref{defn: influenceregion}, and merge the cells accordingly.
	    \State If neighboring influence regions overlap, then merge the cells from both influence regions.
	 \EndFor
	 \State We now have the post-merging grid (at time $t^n$) $\{\tilde{I}_{j_{\text{new}}}(t)\text{ : }j_{\text{new}}=1,2,...,N_{\text{merged}}\}$.
	 
	 \For{$j_{\text{new}}=1,2,\dots,N_{\text{merged}}$}
	    \State $I_{j_{\text{new}}}$ consists of the uniform cells $I_i$ for $i\in\mathcal{S}_{j_{\text{\text{new}}}}$.
		\State Gather the information for $I_{j_{new}}$:
		\State\qquad Take $\nu_{j_{\text{new}}\pm\frac{1}{2}}$ from $\{\nu_{j+\frac{1}{2}}\text{ : }j=0,1,...,N_x\}$.
		\State\qquad Take $\hat{F}_{j_{\text{new}}\pm\frac{1}{2}}(t^n)$ from $\{\hat{F}_{j+\frac{1}{2}}(t^n)\text{ : }j=0,1,...,N_x\}$.
		\State\qquad Compute $\int_{I_{j_{\text{new}}}}{u(x,t^n)dx}\approx\Delta x\sum\limits_{i\in\mathcal{S}_{j_{\text{new}}}}{\bar{u}_i^n}$.
	    \State Compute $\int_{\tilde{I}_{j_{\text{new}}}(t^{(1)})}{u(x,t^{(1)})dx}$ using equation \eqref{eq: SSPRK2_t1_v2}, and hence obtain $\tilde{u}_{j_{\text{new}}}^{(1)}$.
	 \EndFor
	 
	 \For{$j_{\text{new}}=1,2,\dots,N_{\text{merged}}$}
	    \State Compute $\hat{F}_{j_{\text{new}}\pm\frac{1}{2}}(t^{(1)})$ using \textit{nonuniform} ENO3 or WENO-AO3 reconstruction for $u_{j_{\text{new}}+\frac{1}{2}}^{(1),\pm}$.
	    \State Compute $\int_{I_{j_{\text{new}}}^*}{u(x,t^{n+1})dx}$ using equation \eqref{eq: SSPRK2_t2_v2}, and hence obtain $\tilde{u}_{j_{\text{new}}}^{n+1}$.
	 \EndFor

	 \For{$j=1,2,\dots, N_x$}
		\State We have that $I_j\subseteq\bigcup\limits_{j_{\text{new}}=\ell}^{j_{\text{new}}=r}{I_{j_{\text{new}}}^*}$ for some $\ell\leq r$.
		\If{$\ell=r$}
		    \State Compute $\mathlarger{\bar{u}_j^{n+1}\approx\frac{1}{\Delta x}\int\limits_{I_j}{P_{\ell}(x;t^{n+1})dx}}$.
		\Else
		    \State Compute $\mathlarger{\bar{u}_j^{n+1}\approx\frac{1}{\Delta x}\left(\int\limits_{I_j\cap I_{\ell}^*}{P_{\ell}(x;t^{n+1})dx} + \sum\limits_{j_{\text{new}}=\ell+1}^{r-1}{\Delta\tilde{x}_{j_{\text{new}}}^{n+1}\tilde{u}_{j_{\text{new}}}^{n+1}} + \int\limits_{I_j\cap I_{r}^*}{P_{r}(x;t^{n+1})dx}\right)}$.
		 \EndIf
	 \EndFor
	\end{algorithmic}
\end{algorithm}

\begin{rem}
The backward EL-RK-FV algorithm in \cite{nakao2022eulerian} does not use nonuniform WENO-AO schemes for spatial reconstruction over nonuniform grids. Instead, it utilizes a nonuniform-to-uniform transformation that allows one to apply uniform WENO-AO when the velocity field is smooth.
\end{rem}

\begin{rem}
The first- and higher-order algorithms for the forward EL-RK-FV can be easily extended to handle convection-diffusion equations of the form $u_t + \nabla\cdot\mathbf{F}(u) = \epsilon\Delta u + g(\mathbf{x},t)$ by using implicit-explicit Runge-Kutta time discretizations \cite{ascher1997implicit}. Similar to the backward EL-RK-FV scheme in \cite{nakao2022eulerian}, the diffusion term can be discretized using nonuniform linear reconstructions, and the stiff terms can be implicitly updated.
\end{rem}


\subsection{Two-dimensional problems via Strang splitting}

We use dimensional splitting, namely the second-order Strang splitting, to extend the forward EL-RK-FV algorithm to higher dimensions. Although only Burgers' equation is herein addressed for simplicity, we can consider a general two-dimensional hyperbolic conservation law

\begin{equation}\label{eq: scalar2d}
	u_t + f(u)_x + g(u)_y = 0.
\end{equation}

Consider the one-dimensional discretizations in $x$ and $y$ respectively given by
\[0=x_{\frac12}<x_{\frac32}<...<x_{N_x-\frac12}<x_{N_x+\frac12}=2\pi,\qquad 0=y_{\frac12}<y_{\frac32}<...<y_{N_y-\frac12}<y_{N_y+\frac12}=2\pi.\]

The spatial domain $[0,2\pi]\times[0,2\pi]$ is discretized into $N_xN_y$ cells defined by $I_{i,j}\coloneqq I_i\times I_j$ with cell centers $x_{i,j} = ((x_{i-\frac12}+x_{i+\frac12})/2,(y_{j-\frac12}+y_{j+\frac12})/2)$. The cell averages of the solution are defined by
\begin{equation}
	\tilde{\bar{u}}_{i,j}(t) \coloneqq \frac{1}{\Delta x_i\Delta y_i}\iint\limits_{I_{i,j}}{u(x,y,t)dxdy},
\end{equation}
where $\Delta x_i=x_{i+\frac12}-x_{i-\frac12}$ and $\Delta y_j=y_{j+\frac12}-y_{j-\frac12}$. Further let $\tilde{\bar{u}}_{i,j}^n\approx\tilde{\bar{u}}_{i,j}(t^n)$. Only in this subsection will the superscript tilde denote integration (i.e,. averaging) in $y$; it does not denote nonuniform cell averages like in the previous sections. Dimensional splitting methods solve equation \eqref{eq: scalar2d} by alternating between solving the easier one-dimensional problems

\begin{subequations}
\begin{equation}\label{eq: splitting_x}
	u_t + f(u)_x = 0,
\end{equation}
\begin{equation}\label{eq: splitting_y}
	u_t + g(u)_y = 0.
\end{equation}
\end{subequations}

Dimensional splitting in the finite volume framework requires a bit more care than in the finite difference framework since we are evolving cell averages instead of point values. If we are to evolve the two-dimensional cell averages in a dimension-by-dimension manner, then we must integrate with respect to one variable while holding the other variable constant at a point. Define the uniform \textit{interval averages} at time $t$ with one variable fixed by

\begin{subequations}
\begin{equation}\label{eq: intervalaverages_x}
	\bar{u}_{i|y} \coloneqq \psi_i(y,t) = \frac{1}{\Delta x_i}\int\limits_{I_i}{u(x,y,t)dx},
\end{equation}
\begin{equation}\label{eq: intervalaverages_y}
	\tilde{u}_{j|x} \coloneqq \phi_j(x,t) = \frac{1}{\Delta y_j}\int\limits_{I_j}{u(x,y,t)dy}.
\end{equation}
\end{subequations}

Observe that we view the interval average in $x$ as a function of $y$ and $t$, and the interval average in $y$ as a function of $x$ and $t$. The idea of dimensional splitting in the finite volume framework is as follows: fix the solution in $y$ at $L$ Gauss-Legendre quadrature nodes in $I_j$, solve equation \eqref{eq: splitting_x} by evolving the $x$-interval averages \eqref{eq: intervalaverages_x} at each quadrature node, and apply the quadrature by taking a weighted sum of the $x$-interval averages at each Gauss-Legendre node to recover the updated cell average. Repeat in the other dimension, and alternate as needed for the dimensional splitting method. We illustrate going between $x$-interval averages and cell averages.

\begin{align}\label{eq: cellavg2intavg}
\begin{split}
	\tilde{\bar{u}}_{i,j} &= \frac{1}{\Delta y_j}\int\limits_{I_j}{\psi(y)dy}\\
	&= \frac{1}{2}\int\limits_{-1}^{1}{\psi\left(y_{j-\frac12}+\frac{\Delta y}{2}(y'+1)\right)dy'}\\
	&\approx \frac{1}{2}\sum\limits_{l=1}^{L}{w_l\psi(y_l^{(j)})},\\
	&= \frac{1}{2}\sum\limits_{l=1}^{L}{w_l\bar{u}_{i|l}},
\end{split}
\end{align}
where $y_l^{(j)} = y_{j-\frac12}+\frac{\Delta y}{2}(\xi_l+1)$ and $w_l$, $l=1,...,L$ are the Gauss-Legendre quadrature nodes and weights in $I_j$, respectively. Notice that by the first line of equation \eqref{eq: cellavg2intavg}, we can compute $\psi(y)=\bar{u}_{i|y}$ at each quadrature node using ENO or WENO-AO reconstruction polynomials; see Algorithm 3 in \cite{nakao2022eulerian}. These can be used to initialize the scheme. For the sake of brevity, we refer the reader to Algorithm 4 in \cite{nakao2022eulerian} for the Strang splitting algorithm applied to the EL-RK-FV method. The process is identical with the exception that we will use the forward EL-RK-FV scheme instead of the backward EL-RK-FV scheme to evolve the cell averages.


\section{Numerical tests}
\label{sec:numerics}

In this section, we present the numerical performance of the proposed forward EL-RK-FV scheme. The scheme's high-order accuracy is demonstrated on a linear advection equation with a smooth solution, as well as Burgers' equation both before and after shock formations. We also present results showing the scheme's ability to capture shocks, rarefaction waves, and interactions between shocks and rarefaction waves in both the one-dimensional and two-dimensional cases. All numerical tests further demonstrate relatively large allowable time-stepping sizes. We compare the performance of the forward EL-RK-FV scheme using forward Euler method and optimal SSP-RK methods \eqref{eq: SSPRK2_v2} and \eqref{eq: SSPRK3_v2}, which we respectively denote RK1, RK2, and RK3. The time discretization is coupled with either third-order nonuniform ENO or WENO-AO reconstruction, which we respectively denote ENO3 and WENO3.

We note that in our numerical tests, ENO3 and WENO3 perform comparably, with ENO3 perfoming marginally better on some occasions. This is because our third-order nonuniform WENO-AO reconstruction (as shown in the appendix) borrows the idea from \cite{wenoao2016} in which a larger center stencil is used in combination with smaller left/right stencils. In our case, we use a centered 3-stencil $\{I_{j-1},I_j,I_{j+1}\}$ and left/right-biased 2-stencils $\{I_{j-1},I_j\}$ and $\{I_j,I_{j+1}\}$. Whereas, ENO3 uses 3-stencils for the centered and left/right-biased stencils, thus considering more information than our WENO3.

Strang splitting is used for dimensional splitting, but higher-order splitting methods \cite{forest1990fourth,yoshida1990construction} can also be used. Lastly, the time-stepping sizes in one and two dimensions are respectively defined as

\begin{subequations}
\begin{equation}\label{eq: CFL1d}
	\Delta t = \frac{\text{CFL}\cdot\Delta x}{\max|f'(u)|},
\end{equation}
\\
\begin{equation}\label{eq: CFL2d}
	\Delta t = \frac{\text{CFL}}{\mathlarger{\frac{\max|f'(u)|}{\Delta x} + \frac{\max|g'(u)|}{\Delta y}}},
\end{equation}
\end{subequations}
\ \\
where ``CFL" is defined to compare with the classical CFL number from an Eulerian approach. Recall from Remark \ref{rem: timesteprestriction} that in the presence of troubled cells, the time-stepping size is restricted by $\Delta t<\lambda\Delta x$ according to equation \eqref{eq: timeconstraint}. Although the stability proof in \cite{yang2023stability} requires a strict inequality, it is possible that equality is allowable. Some of our numerical experiments demonstrated stability if $\Delta t\leq\lambda\Delta x$. We note that in the two-dimensional case, this restriction is

\begin{equation}\label{eq: timeconstraint2d}
	\Delta t<\frac{4\min{\{\Delta x,\Delta y\}}}{\max{u_0(x,y)}-\min{u_0(x,y)}}.
\end{equation}
\\
\begin{exa}
(1D transport with variable coefficient)
\begin{equation}\label{eq: varcoeff}
	u_t + (\sin{(x)}u)_x = 0,\qquad x\in[0,2\pi],
\end{equation}
with periodic boundary conditions and exact solution $u(x,t) = \sin{(2\arctan{(e^{-t}\tan{(x/2)})})}/\sin{(x)}$. This example is meant to verify the effectiveness and order of convergence for the proposed high-order forward EL-RK-FV scheme. Since shocks do not form for equation \eqref{eq: varcoeff}, the conditions in Definitions \ref{defn: troubledcell} and \ref{defn: influenceregion} do not activate and we are not restricted by inequality \eqref{eq: timeconstraint}. The solution is simply evolved over each space-time region $\Omega_j$ using SSP-RK methods and nonuniform reconstruction. Referring to Table \ref{tab: varcoeff_errortable}, we observe the expected high-order convergence when solving up to final time $T_f=1$ with a relatively large time-stepping size using $\text{CFL}=3.2$. Due to the large CFL number, the temporal error dominates the spatial error. This is further seen in Figure \ref{fig: varcoeff_errorplot}, where we verify the high-order accuracy by fixing the mesh $N_x=200$, solving up to final time $T_f=0.5$, and varying the CFL number from 0.01 to 10. Note that the temporal error starts to dominate at $\text{CFL}=0.2$ and $\text{CFL}=1$ when using RK2 and RK3, respectively. Furthermore, the $L^1$ errors are nearly indistinguishable for ENO3 and WENO3.

\begin{table}[t!] 
\begin{center}
\caption{Convergence study with spatial mesh refinement for equation \eqref{eq: varcoeff} with initial condition $u_0(x)=1$, solving up to $T_f=1$ with $\Delta t = 3.2\Delta x$ ($\text{CFL}=3.2$).}
\label{tab: varcoeff_errortable}
\begin{tabular}{|*{7}{c|}} 
\hline 
\multicolumn{7}{|c|}{WENO3 Reconstruction}\\
\hline
\multicolumn{1}{|c|}{}&\multicolumn{2}{|c|}{RK1}&\multicolumn{2}{|c|}{RK2}&\multicolumn{2}{|c|}{RK3}\\
\hline
$N_x$&$L^1$ Error&Order&$L^1$ Error&Order&$L^1$ Error&Order\\
\hline
100&6.49E-02&-&3.53E-03&-&2.11E-04&-\\
200&3.17E-02&1.00&7.95E-04&2.15&2.45E-05&3.11\\
300&2.10E-02&1.01&3.43E-04&2.07&7.16E-06&3.03\\
400&1.57E-02&1.01&1.90E-04&2.05&3.00E-06&3.02\\
\hline
\multicolumn{7}{|c|}{ENO3 Reconstruction}\\
\hline
\multicolumn{1}{|c|}{}&\multicolumn{2}{|c|}{RK1}&\multicolumn{2}{|c|}{RK2}&\multicolumn{2}{|c|}{RK3}\\
\hline
$N_x$&$L^1$ Error&Order&$L^1$ Error&Order&$L^1$ Error&Order\\
\hline
100&6.49E-02&-&3.53E-03&-&2.08E-04&-\\
200&3.17E-02&1.00&7.95E-04&2.15&2.45E-05&3.08\\
300&2.10E-02&1.01&3.43E-04&2.07&7.17E-06&3.03\\
400&1.57E-02&1.01&1.90E-04&2.05&3.01E-06&3.02\\
\hline
\end{tabular} 
\end{center} 
\end{table} 

\begin{figure}[t!]
	\centering
	\includegraphics[width=0.65\textwidth]{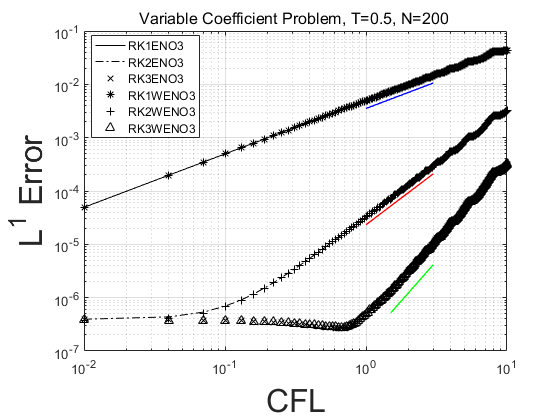}
	\caption{Equation \eqref{eq: varcoeff} with initial condition $u_0(x)=1$, solving up to final time $T_f=0.5$ with mesh $N_x=200$.}
	\label{fig: varcoeff_errorplot}
\end{figure}
\end{exa}

\begin{exa}
(1D Burgers' equation with smooth initial data)
\begin{equation}\label{eq: burgers}
	u_t + \left(\frac{u^2}{2}\right)_x = 0,\qquad x\in[0,2\pi],
\end{equation}
with periodic boundary conditions and initial condition $u_0(x) = \sin{(x)}$. We test the order of convergence both before and after the breaking time $t_B=1$ when the shock forms. The exact solution from the method of characteristics is computed using Newton's method. Table \ref{tab: burgers_errortable_beforeshock} shows the expected high-order convergence when solving up to final time $T_f=0.5$ using $\text{CFL}=3.2$. Figure \ref{fig: burgers_errorplot_beforeshock} shows the change in the $L^1$ error as the CFL number varies from 0.01 to 30 with a fixed mesh $N_x=200$. As with equation \eqref{eq: varcoeff}, the time-stepping size is not restricted by inequality \eqref{eq: timeconstraint} since the shock has not yet formed. As seen in Table \ref{tab: burgers_errortable_beforeshock} and Figure \ref{fig: burgers_errorplot_beforeshock}, we observe that the forward EL-RK-FV scheme performs slightly better when using ENO3 over WENO3. Also, the $L^1$ error is the same regardless of the RK method used since the space-time partition lines defined by equation \eqref{eq: forward_characteristics} exactly trace the characteristics for Burgers' equation. Hence, $F_{j_{\text{new}}+\frac{1}{2}}(t)=0$ and we do not expect to see any temporal error.

\begin{table}[t!] 
\begin{center}
\caption{Convergence study with spatial mesh refinement for equation \eqref{eq: burgers} with initial condition $u_0(x)=\sin{(x)}$, solving up to $T_f=0.5$ with $\Delta t = 3.2\Delta x$ ($\text{CFL}=3.2$).}
\label{tab: burgers_errortable_beforeshock}
\begin{tabular}{|*{7}{c|}} 
\hline 
\multicolumn{7}{|c|}{WENO3 Reconstruction}\\
\hline
\multicolumn{1}{|c|}{}&\multicolumn{2}{|c|}{RK1}&\multicolumn{2}{|c|}{RK2}&\multicolumn{2}{|c|}{RK3}\\
\hline
$N_x$&$L^1$ Error&Order&$L^1$ Error&Order&$L^1$ Error&Order\\
\hline
100&4.00E-05&-&4.00E-05&-&4.00E-05&-\\
200&6.23E-06&2.68&6.23E-06&2.68&6.23E-06&2.68\\
300&2.39E-06&2.36&2.39E-06&2.36&2.39E-06&2.36\\
400&5.31E-07&5.24&5.31E-06&5.24&5.31E-07&5.24\\
\hline
\multicolumn{7}{|c|}{ENO3 Reconstruction}\\
\hline
\multicolumn{1}{|c|}{}&\multicolumn{2}{|c|}{RK1}&\multicolumn{2}{|c|}{RK2}&\multicolumn{2}{|c|}{RK3}\\
\hline
$N_x$&$L^1$ Error&Order&$L^1$ Error&Order&$L^1$ Error&Order\\
\hline
100&4.72E-06&-&4.72E-06&-&4.72E-06&-\\
200&6.27E-07&2.91&6.27E-07&2.91&6.27E-07&2.91\\
300&2.00E-07&2.82&2.00E-07&2.82&2.00E-07&2.82\\
400&8.47E-08&2.98&8.47E-08&2.98&8.47E-08&2.98\\
\hline
\end{tabular} 
\end{center} 
\end{table}

\begin{figure}[t!]
	\centering
	\includegraphics[width=0.65\textwidth]{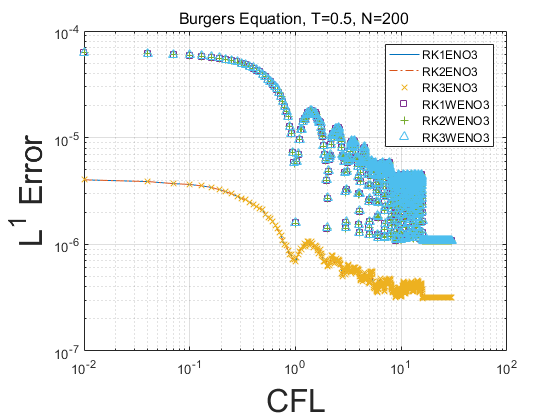}
	\caption{Equation \eqref{eq: burgers} with initial condition $u_0(x)=\sin{(x)}$, solving up to final time $T_f=0.5$ with mesh $N_x=200$.}
	\label{fig: burgers_errorplot_beforeshock}
\end{figure}

When testing for convergence after the shock has formed, the time-stepping size must satisfy the constraint given in \eqref{eq: timeconstraint}, which for this problem is $\Delta t< 2\Delta x$. Table \ref{tab: burgers_errortable_aftershock} shows the expected high-order convergence when solving up to final time $T_f=1.3$ using $\text{CFL}=1.95$. Due to the presence of the shock at $x=\pi$, we only compute the $L^1$ error over the region $[0,\pi-0.1]\cup[\pi+0.1,2\pi]$. Figure \ref{fig: burgers_errorplot_aftershock} shows the change in the $L^1$ error as the CFL number varies from 0.1 to 3.5. Although the scheme appears stable for CFL numbers slightly larger than 2, we lose stability by $\text{CFL}=3$. As with the pre-shock test, we observe that the scheme performs slightly better with ENO3 over WENO3, and the time discretization does not contribute any error since the characteristics are traced exactly.

\begin{table}[t!] 
\begin{center}
\caption{Convergence study with spatial mesh refinement for equation \eqref{eq: burgers} with initial condition $u_0(x)=\sin{(x)}$, solving up to $T_f=1.3$ with $\Delta t = 1.95\Delta x$ ($\text{CFL}=1.95$).}
\label{tab: burgers_errortable_aftershock}
\begin{tabular}{|*{7}{c|}} 
\hline 
\multicolumn{7}{|c|}{WENO3 Reconstruction}\\
\hline
\multicolumn{1}{|c|}{}&\multicolumn{2}{|c|}{RK1}&\multicolumn{2}{|c|}{RK2}&\multicolumn{2}{|c|}{RK3}\\
\hline
$N_x$&$L^1$ Error&Order&$L^1$ Error&Order&$L^1$ Error&Order\\
\hline
100&3.43E-05&-&2.93E-04&-&3.48E-04&-\\
200&8.27E-06&2.05&1.08E-05&4.76&1.05E-05&5.05\\
300&2.62E-06&2.83&2.63E-06&3.49&2.62E-06&3.41\\
400&5.26E-07&5.58&5.25E-07&5.58&5.25E-07&5.58\\
\hline
\multicolumn{7}{|c|}{ENO3 Reconstruction}\\
\hline
\multicolumn{1}{|c|}{}&\multicolumn{2}{|c|}{RK1}&\multicolumn{2}{|c|}{RK2}&\multicolumn{2}{|c|}{RK3}\\
\hline
$N_x$&$L^1$ Error&Order&$L^1$ Error&Order&$L^1$ Error&Order\\
\hline
100&4.07E-06&-&2.59E-04&-&3.08E-04&-\\
200&6.35E-07&2.68&9.25E-07&8.13&1.33E-06&7.86\\
300&1.74E-07&3.19&1.74E-07&4.11&1.74E-07&5.01\\
400&7.00E-08&3.17&7.00E-08&3.17&7.00E-08&3.17\\
\hline
\end{tabular} 
\end{center} 
\end{table} 
\begin{figure}[t!]
	\centering
	\includegraphics[width=0.65\textwidth]{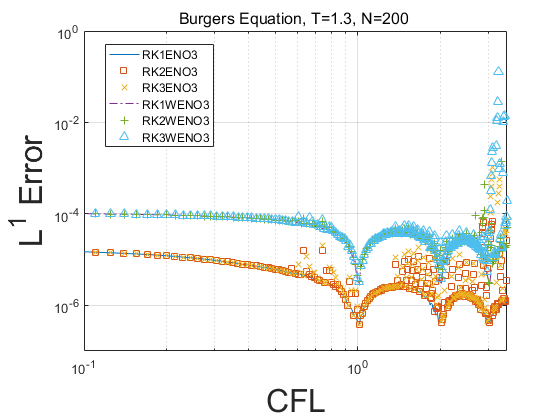}
	\caption{Equation \eqref{eq: burgers} with initial condition $u_0(x)=\sin{(x)}$, solving up to final time $T_f=1.3$ with mesh $N_x=200$.}
	\label{fig: burgers_errorplot_aftershock}
\end{figure}
\end{exa}

\begin{exa}
(1D Burgers' equation with smooth initial data: shock-rarefaction interaction)\\
\ \\
We test Burgers' equation \eqref{eq: burgers} with smooth initial condition $u_0(x) = 1+2\sin{(x)}$ and periodic boundary conditions. A shock forms at the breaking time $t_B=0.5$, so we solve up to final time $T_f=1.3$ with mesh $N_x=100$ and time-stepping size $\Delta t = 0.99\Delta x$. Due to the shock, we must satisfy the stability constraint $\Delta t<\Delta x$ ($\text{CFL}=3$) from inequality \eqref{eq: timeconstraint}. Figure \ref{fig: burgers1d_sinshock} presents the solution plots comparing different RK schemes with ENO3 and WENO3. Regardless of the RK scheme used, ENO3 and WENO3 give comparable results. RK1 (forward Euler method) does a slightly better job capturing the shock since, as a one-step method, it requires fewer transition points for the shocks. Further note that a rarefaction wave forms at the boundary. By time $t=1.3$, the shock and rarefaction wave have interacted, with the rarefaction wave slowly consuming the shock as the solution evolves.

\begin{figure}[t!]
\centering
\begin{minipage}[b]{0.32\linewidth}
	\centering
	\includegraphics[width=\textwidth]{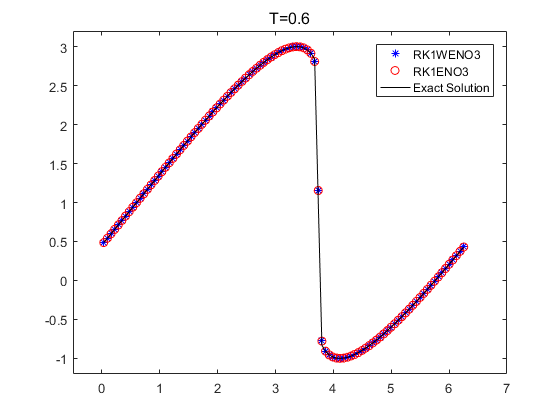}\\
	\includegraphics[width=\textwidth]{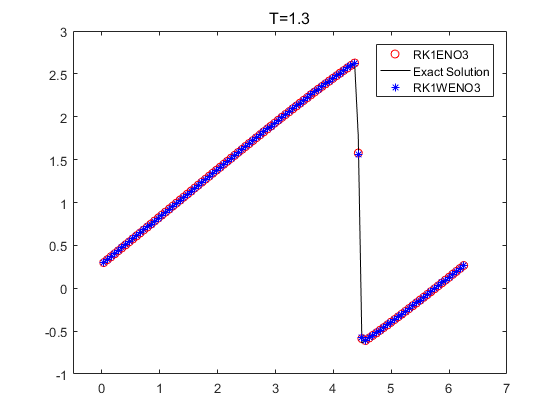}
\end{minipage}
\begin{minipage}[b]{0.32\linewidth}
	\centering
	\includegraphics[width=\textwidth]{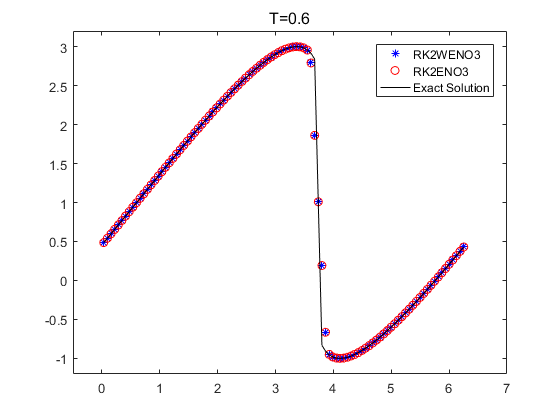}\\
	\includegraphics[width=\textwidth]{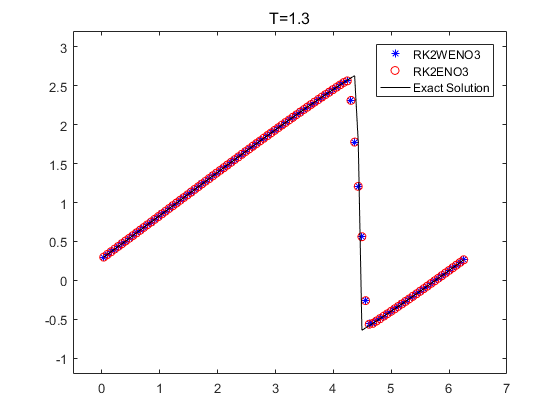}
\end{minipage}
\begin{minipage}[b]{0.32\linewidth}
	\centering
	\includegraphics[width=\textwidth]{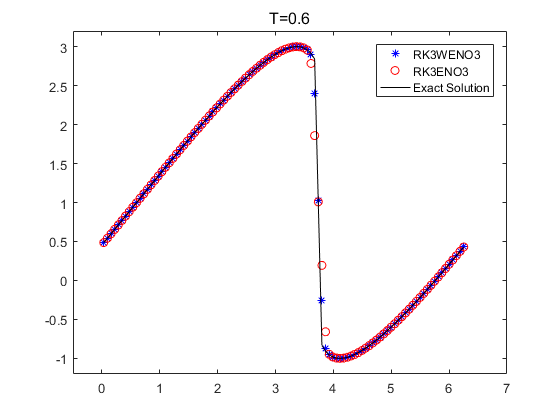}\\
	\includegraphics[width=\textwidth]{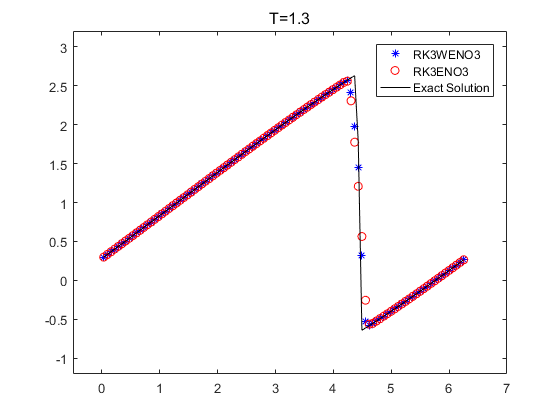}
\end{minipage}
\caption{Equation \eqref{eq: burgers} with initial condition $u_0(x) = 1+2\sin{(x)}$. Mesh $N_x=100$, time-stepping size $\Delta t = 0.99\Delta x$ ($\text{CFL}=2.97$), and plotting the solution at times $t=0.6$ and $t=1.3$. (Top row) $t=0.6$, (bottom row) $t=1.3$.}
\label{fig: burgers1d_sinshock}
\end{figure}
\end{exa}

\begin{exa}
(1D Burgers' equation with discontinuous initial data: shock)\\
\ \\
We test the Riemann problem with Burgers' equation \eqref{eq: burgers} using initial condition
\begin{equation}\label{eq: 1DRiemann_shock}
	u_0(x) = 
	\begin{cases}
		4,&x\leq 0,\\
		0,&x>0.
	\end{cases}
\end{equation}
Solving over the domain $[-\pi,\pi]$ and assuming Dirichlet boundary conditions, a shock instantaneously forms at $x=0$ and travels to the right with speed $\nu=2$. We solve up to final time $T_f=1.2$ with mesh $N_x=100$ and time-stepping size $\Delta t = \Delta x$ ($\text{CFL}=4$). Although constraint \eqref{eq: timeconstraint} requires $\Delta t<\Delta x$, we observe stability if $\Delta t\leq\Delta x$. Figure \ref{fig: burgers1d_shock} presents the solution plots comparing different RK schemes with ENO3 and WENO3. Regardless of the RK scheme used, ENO3 and WENO3 give comparable results. However, all the schemes do a good job capturing the shock despite the large CFL number. Similar to the previous example, RK1 (forward Euler method) does a slightly better job capturing the shock since, as a one-step method, it requires fewer transition points for the shocks.

\begin{figure}[t!]
\centering
\begin{minipage}[b]{0.32\linewidth}
	\centering
	\includegraphics[width=\textwidth]{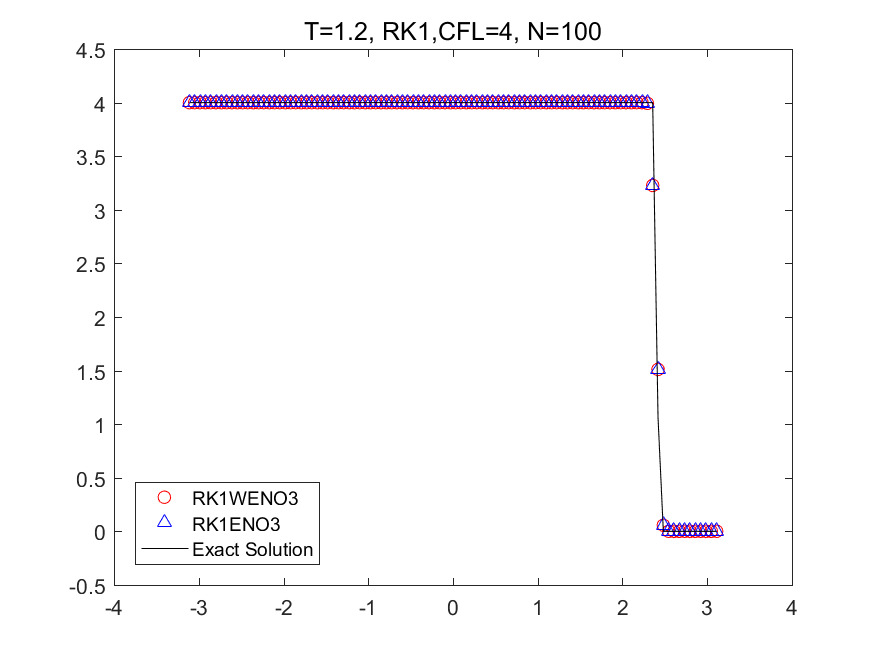}
\end{minipage}
\begin{minipage}[b]{0.32\linewidth}
	\centering
	\includegraphics[width=\textwidth]{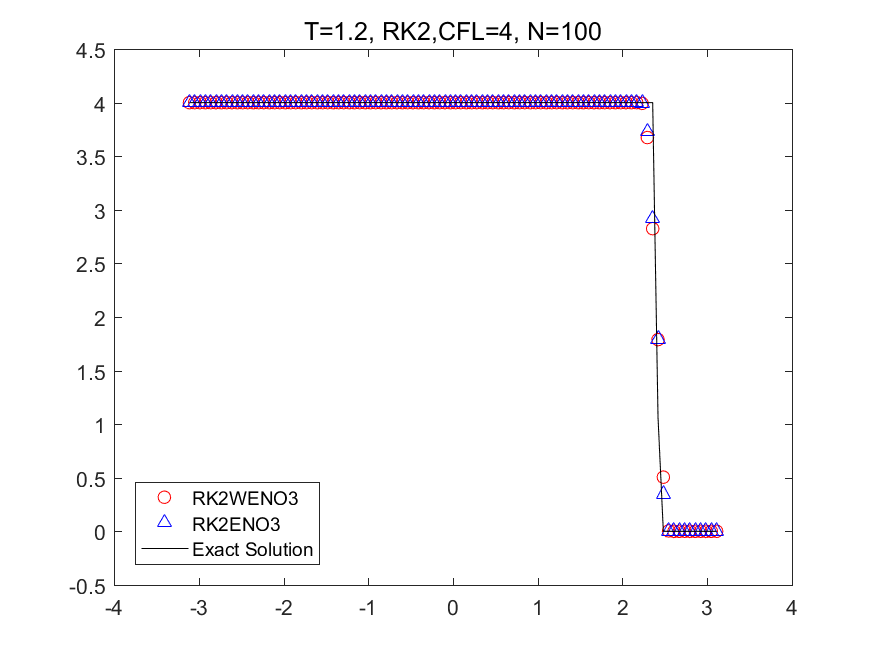}
\end{minipage}
\begin{minipage}[b]{0.32\linewidth}
	\centering
	\includegraphics[width=\textwidth]{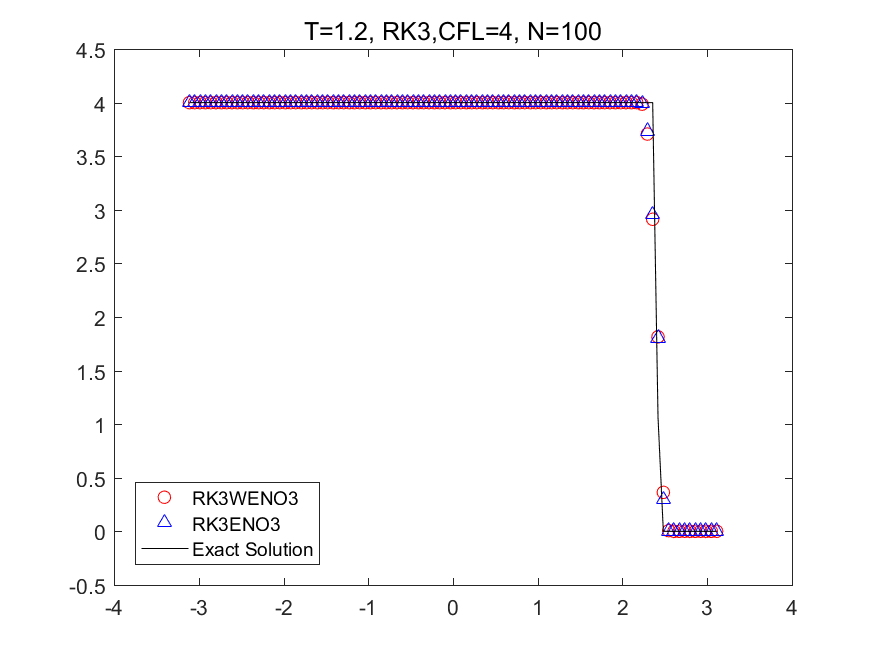}
\end{minipage}
\caption{Equation \eqref{eq: burgers} with initial condition \eqref{eq: 1DRiemann_shock}, solving up to time $T_f=1.2$ with mesh $N_x=100$ and time-stepping size $\Delta t = \Delta x$ ($\text{CFL}=4$). (Left) RK1, (middle) RK2, (right) RK3.}
\label{fig: burgers1d_shock}
\end{figure}
\end{exa}

\begin{exa}
(1D Burgers' equation with discontinuous initial data: rarefaction wave)\\
\ \\
We test the Riemann problem with Burgers' equation \eqref{eq: burgers} using initial condition
\begin{equation}\label{eq: 1DRiemann_rarefaction}
	u_0(x) = 
	\begin{cases}
		-2,&x\leq 0,\\
		2,&x>0.
	\end{cases}
\end{equation}
Solving over the domain $[-\pi,\pi]$ and assuming Dirichlet boundary conditions, a rarefaction wave forms at $x=0$ and travels in each direction with equal speed. We solve up to final time $T_f=1.2$ with mesh $N_x=100$ and time-stepping size $\Delta t=1.6\Delta x$ ($\text{CFL}=3.2$). Since a single rarefaction wave forms and does not trigger the merging procedure, we can use time-stepping sizes that exceed constraint \eqref{eq: timeconstraint} which says $\Delta t<\Delta x$ ($\text{CFL}=2$). As seen in Figure \ref{fig: burgers1d_rarefaction}, all the schemes evolve the solution well despite the large CFL number, with RK3 performing the best. As with the single shock case in initial condition \eqref{eq: 1DRiemann_shock}, ENO3 and WENO3 produce comparable results.

\begin{figure}[t!]
\centering
\begin{minipage}[b]{0.32\linewidth}
	\centering
	\includegraphics[width=\textwidth]{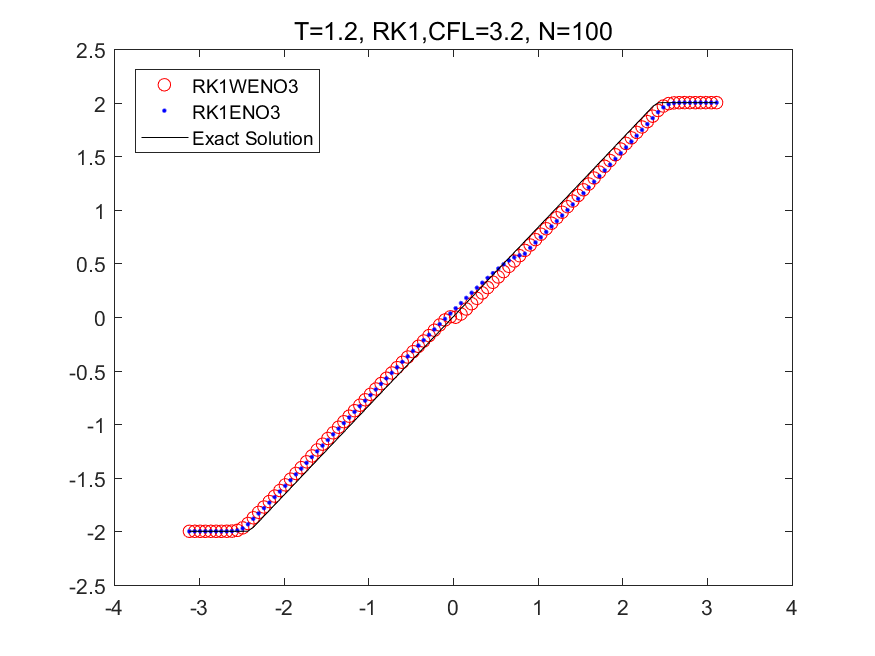}
\end{minipage}
\begin{minipage}[b]{0.32\linewidth}
	\centering
	\includegraphics[width=\textwidth]{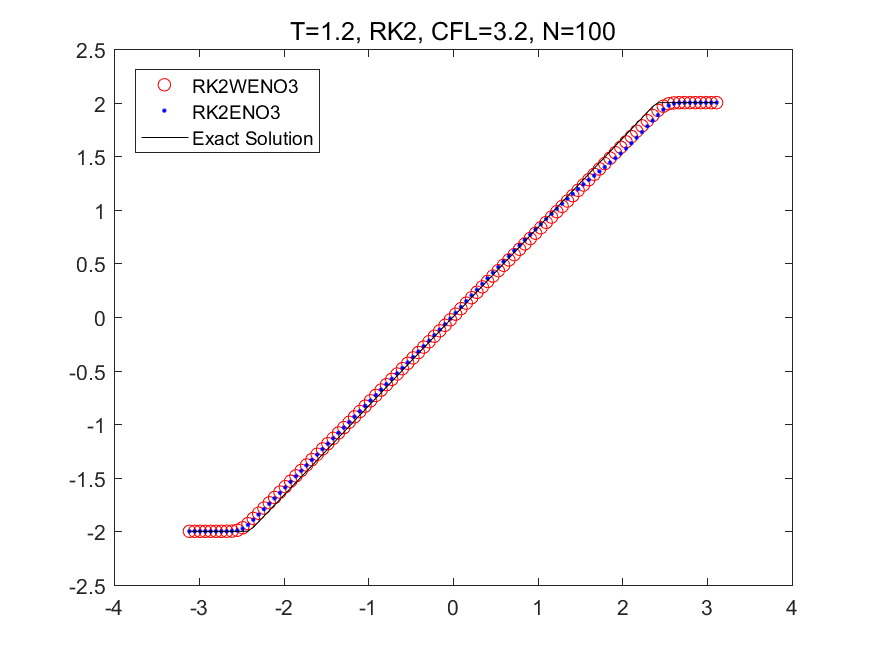}
\end{minipage}
\begin{minipage}[b]{0.32\linewidth}
	\centering
	\includegraphics[width=\textwidth]{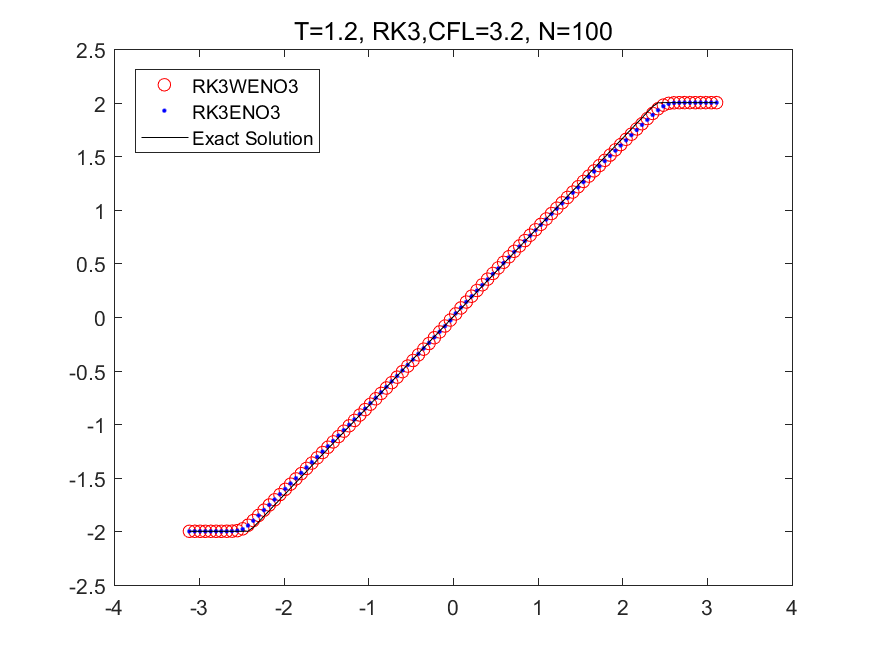}
\end{minipage}
\caption{Equation \eqref{eq: burgers} with initial condition \eqref{eq: 1DRiemann_rarefaction}, solving up to time $T_f=1.2$ with mesh $N_x=100$ and time-stepping size $\Delta t = 1.6\Delta x$ ($\text{CFL}=3.2$). (Left) RK1, (middle) RK2, (right) RK3.}
\label{fig: burgers1d_rarefaction}
\end{figure}
\end{exa}

\begin{exa}
(2D Burgers' equation with continuous initial data)
\begin{equation}\label{eq: burgers2d}
	u_t + \left(\frac{u^2}{2}\right)_x + \left(\frac{u^2}{2}\right)_y = 0,\qquad (x,y)\in[0,2]\times[0,2],
\end{equation}
with homogeneous Dirichlet boundary conditions and initial condition
\begin{equation}\label{eq: burgers2d_sinIC}
	u_0(x,y) = 
	\begin{cases}
		\sin^2{(\pi x)}\sin^2{(\pi y)},&(x,y)\in[0,1]\times[0,1],\\
		0,&\text{otherwise}.
	\end{cases}
\end{equation}
The solution is plotted at times $t=1$, $t=2$, and $t=3$, with mesh $N_x=N_y=100$ and time-stepping size $\Delta t=3.9\Delta x$ ($\text{CFL}=7.8$). Since a shock is generated, we must satisfy constraint \eqref{eq: timeconstraint2d}, which for this problem is $\Delta t<4\Delta x$ ($\text{CFL}=8$). We use Strang splitting for dimensional splitting, RK3 for time discretization, and ENO3 for spatial reconstruction since ENO3 and WENO3 had comparable results in the one-dimensional tests with ENO3 having slightly smaller error. As seen in Figure \ref{fig: burgers2d_shocksin}, the forward EL-RK-FV scheme captures the discontinuities formed by the shock while using a very large time-stepping size.

\begin{figure}[t!]
\centering
\begin{minipage}[b]{0.32\linewidth}
	\centering
	\includegraphics[width=\textwidth]{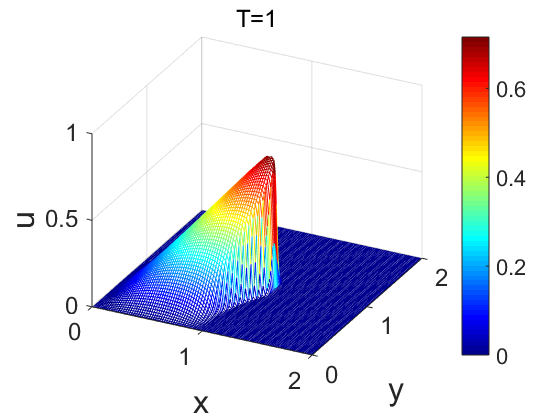}\\
	\includegraphics[width=\textwidth]{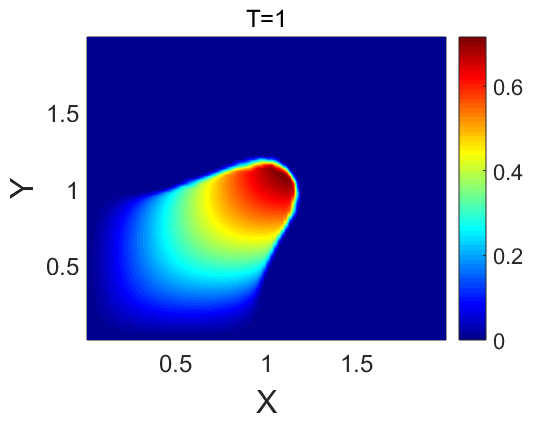}
\end{minipage}
\begin{minipage}[b]{0.32\linewidth}
	\centering
	\includegraphics[width=\textwidth]{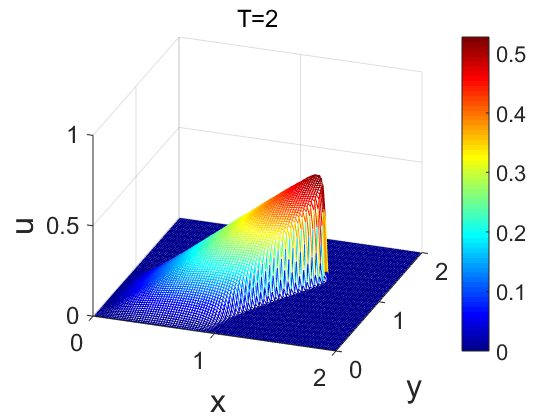}\\
	\includegraphics[width=\textwidth]{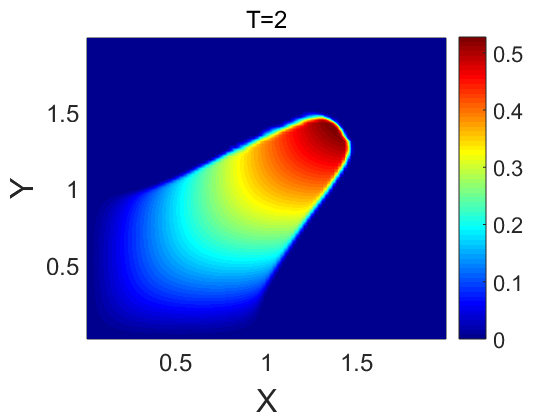}
\end{minipage}
\begin{minipage}[b]{0.32\linewidth}
	\centering
	\includegraphics[width=\textwidth]{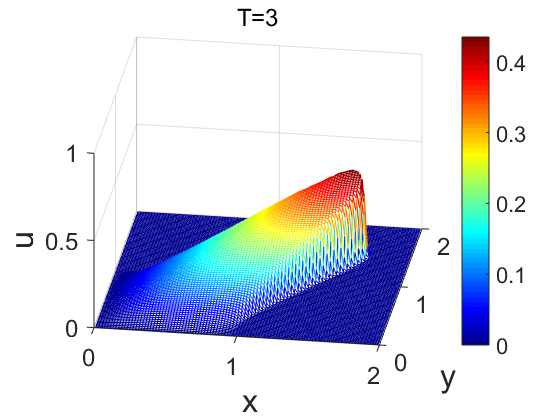}\\
	\includegraphics[width=\textwidth]{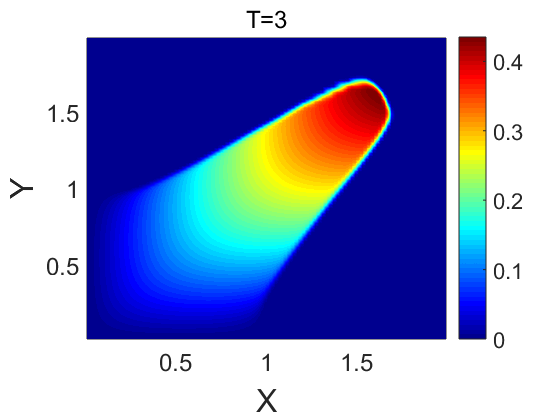}
\end{minipage}
\caption{Equation \eqref{eq: burgers2d} with initial condition \eqref{eq: burgers2d_sinIC}. Mesh $N_x=N_y=100$, time-stepping size $\Delta t = 3.9\Delta x$ ($\text{CFL}=7.8$), RK3 and ENO3, plotting the solution at times $t=1$, $t=2$, and $t=3$. (Top row) mesh plots, (bottom row) contour plots.}
\label{fig: burgers2d_shocksin}
\end{figure}
\end{exa}

\begin{exa}
(2D Burgers' equation with discontinuous initial data: 4 shocks)\\
\ \\
We test the Riemann problem with Burgers' equation \eqref{eq: burgers2d} using initial condition \cite{yoon2008two}
\begin{equation}\label{eq: 2DRiemann_4shock}
	u_0(x,y) = 
	\begin{cases}
		1,&(x,y)\in(0,0.5]\times(0,0.5],\\
		2,&(x,y)\in[-0.5,0]\times[0,0.5],\\
		4,&(x,y)\in[-0.5,0)\times[-0.5,0],\\
		3,&(x,y)\in(0,0.5]\times[-0.5,0),
	\end{cases}
\end{equation}
where the initial values $(a,b,c,d)$ are assigned to each of the four quadrants: $a$ to the first quadrants, $b$ to the second quadrant, $c$ to the third quadrant, and $d$ to the fourth quadrant. (We will use this notation for the remaining Riemann problem examples in this section.) Solving over the domain $[-0.5,0.5]\times[-0.5,0.5]$ and assuming Dirichlet boundary conditions, this initial condition generates four shocks on the major axes at angles $\theta=0$, $\theta=\pi/2$, $\theta=\pi$, and $\theta=3\pi/2$. The solution can be found in \cite{yoon2008two}. Figure \ref{fig: burgers2d_4shock} plots the solution at final time $T_f=0.1$ using mesh $N_x=N_y=100$, and time-stepping size $\Delta t = 1.3\Delta x$ ($\text{CFL}=10.4$). The first-order scheme is compared against the third-order scheme using ENO3 and RK3.

As seen in Figure \ref{fig: burgers2d_4shock}, the $\theta_0$ shock interacts with the $\theta_{3\pi/2}$ shock at the top-left corner of the yellow bottom-right rectangle, and the $\theta_{\pi/2}$ shock interacts with the $\theta_{\pi}$ shock at the bottom-right corner of the light blue top-left rectangle. A new shock is formed at each of these shock-shock interaction points. These two newly formed shocks further interact along the top-right edge of the red lower-left `rectangle'. Observe that the first-order scheme distorts the solution at this top-right edge of the red lower-left `rectangle'. Whereas, the third-order scheme produces no distortion and accurately captures all the discontinuities formed by the shocks.

\begin{figure}[t!]
\centering
\begin{minipage}[b]{0.48\linewidth}
	\centering
	\includegraphics[width=\textwidth]{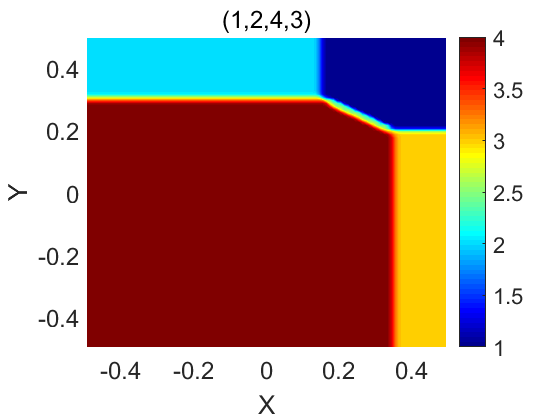}
\end{minipage}
\begin{minipage}[b]{0.48\linewidth}
	\centering
	\includegraphics[width=\textwidth]{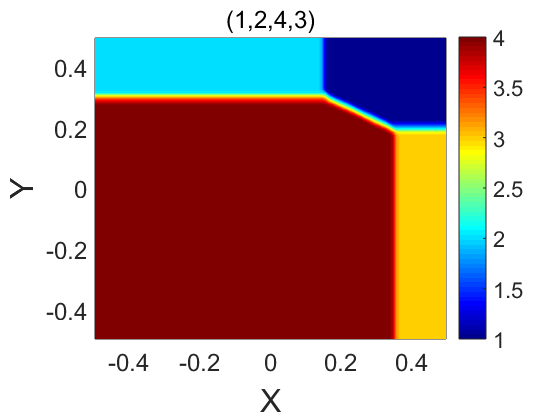}
\end{minipage}
\caption{Equation \eqref{eq: burgers2d} with initial condition \eqref{eq: 2DRiemann_4shock}: 4 shocks. Mesh $N_x=N_y=100$, time-stepping size $\Delta t = 1.3\Delta x$ ($\text{CFL}=10.4$), final time $T_f=0.1$. (Left) first-order scheme, (right) third-order scheme with ENO3 and RK3.}
\label{fig: burgers2d_4shock}
\end{figure}
\end{exa}

\begin{exa}
(2D Burgers' equation with discontinuous initial data: 4 rarefaction waves)\\
\ \\
We test the Riemann problem with Burgers' equation \eqref{eq: burgers2d} using initial condition \cite{yoon2008two}
\begin{equation}\label{eq: 2DRiemann_4rare}
	u_0(x,y) = 
	\begin{cases}
		4,&(x,y)\in(0,0.5]\times(0,0.5],\\
		2,&(x,y)\in[-0.5,0]\times[0,0.5],\\
		1,&(x,y)\in[-0.5,0)\times[-0.5,0],\\
		3,&(x,y)\in(0,0.5]\times[-0.5,0).
	\end{cases}
\end{equation}
Solving over the domain $[-0.5,0.5]\times[-0.5,0.5]$ and assuming Dirichlet boundary conditions, this initial condition generates four rarefaction waves on the major axes at angles $\theta=0$, $\theta=\pi/2$, $\theta=\pi$, and $\theta=3\pi/2$. The solution can be found in \cite{yoon2008two}. Figure \ref{fig: burgers2d_4rare} plots the solution at final time $T_f=0.1$ using mesh $N_x=N_y=100$, and time-stepping size $\Delta t = 1.3\Delta x$ ($\text{CFL}=10.4$). The first-order scheme is compared against the third-order scheme using ENO3 and RK3. Figure \ref{fig: burgers2d_4rare} presents contour plots of the solution as the rarefaction waves evolve and interact with each other. Although both schemes do a job good capturing the solution, the third-order scheme captures the rarefaction waves more sharply and accurately than the first-order scheme.

\begin{figure}[t!]
\centering
\begin{minipage}[b]{0.48\linewidth}
	\centering
	\includegraphics[width=\textwidth]{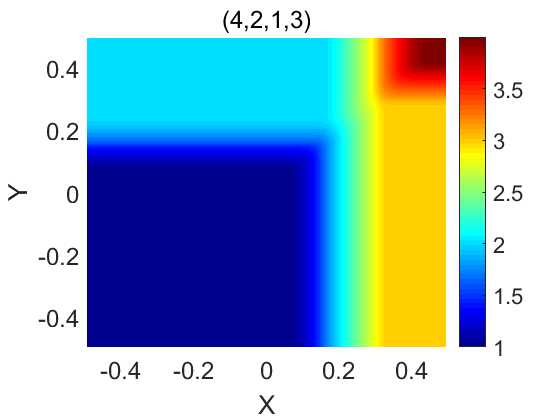}
\end{minipage}
\begin{minipage}[b]{0.48\linewidth}
	\centering
	\includegraphics[width=\textwidth]{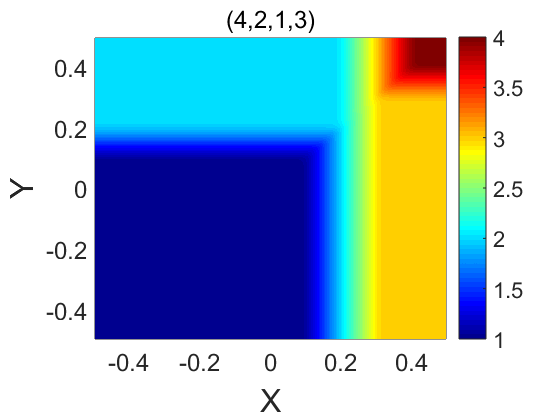}
\end{minipage}
\caption{Equation \eqref{eq: burgers2d} with initial condition \eqref{eq: 2DRiemann_4rare}: 4 rarefaction waves. Mesh $N_x=N_y=100$, time-stepping size $\Delta t = 1.3\Delta x$ ($\text{CFL}=10.4$), final time $T_f=0.1$. (Left) first-order scheme, (right) third-order scheme with ENO3 and RK3.}
\label{fig: burgers2d_4rare}
\end{figure}
\end{exa}

\begin{exa}
(2D Burgers' equation with discontinuous initial data: 1 shock and 3 rarefaction waves)\\
\ \\
We test the Riemann problem with Burgers' equation \eqref{eq: burgers2d} using initial condition \cite{yoon2008two}
\begin{equation}\label{eq: 2DRiemann_1shock3rare}
	u_0(x,y) = 
	\begin{cases}
		4,&(x,y)\in(0,0.5]\times(0,0.5],\\
		3,&(x,y)\in[-0.5,0]\times[0,0.5],\\
		2,&(x,y)\in[-0.5,0)\times[-0.5,0],\\
		1,&(x,y)\in(0,0.5]\times[-0.5,0).
	\end{cases}
\end{equation}
Solving over the domain $[-0.5,0.5]\times[-0.5,0.5]$ and assuming Dirichlet boundary conditions, this initial condition generates one shock at angle $\theta=3\pi/2$, and three rarefaction waves at angles $\theta=0$, $\theta=\pi/2$, and $\theta=\pi$. The solution can be found in \cite{yoon2008two}. Figure \ref{fig: burgers2d_1shock3rare} plots the solution at final time $T_f=0.1$ using mesh $N_x=N_y=100$, and time-stepping size $\Delta t = 1.3\Delta x$ ($\text{CFL}=10.4$). The first-order scheme is compared against the third-order scheme using ENO3 and RK3.

As seen in Figure \ref{fig: burgers2d_1shock3rare}, as the $\theta_{3\pi/2}$ shock travels to the right, the $\theta_0$ rarefaction wave travels up while simultaneously interacting with the shock; this occurs on the upper edge of the blue lower-right `rectangle'. Figure \ref{fig: burgers2d_1shock3rare} also shows how the other rarefaction waves interact as the solution evolves. As with the 4 rarefaction waves example, the third-order scheme captures the rarefaction waves more sharply and accurately than the first-order scheme.

\begin{figure}[t!]
\centering
\begin{minipage}[b]{0.48\linewidth}
	\centering
	\includegraphics[width=\textwidth]{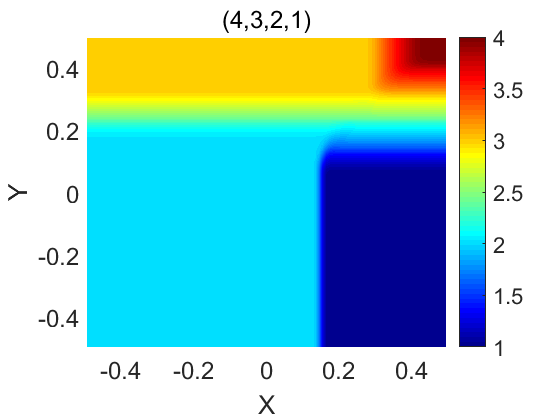}
\end{minipage}
\begin{minipage}[b]{0.48\linewidth}
	\centering
	\includegraphics[width=\textwidth]{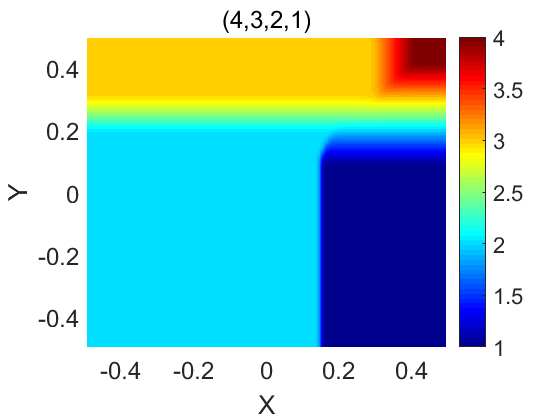}
\end{minipage}
\caption{Equation \eqref{eq: burgers2d} with initial condition \eqref{eq: 2DRiemann_1shock3rare}: 1 shock and 3 rarefaction waves. Mesh $N_x=N_y=100$, time-stepping size $\Delta t = 1.3\Delta x$ ($\text{CFL}=10.4$), final time $T_f=0.1$. (Left) first-order scheme, (right) third-order scheme with ENO3 and RK3.}
\label{fig: burgers2d_1shock3rare}
\end{figure}
\end{exa}

\begin{exa}
(2D Burgers' equation with discontinuous initial data: 2 shocks and 2 rarefaction waves)\\
\ \\
We test the Riemann problem with Burgers' equation \eqref{eq: burgers2d} using initial condition \cite{yoon2008two}
\begin{equation}\label{eq: 2DRiemann_2shocks2rare}
	u_0(x,y) = 
	\begin{cases}
		1,&(x,y)\in(0,0.5]\times(0,0.5],\\
		3,&(x,y)\in[-0.5,0]\times[0,0.5],\\
		2,&(x,y)\in[-0.5,0)\times[-0.5,0],\\
		4,&(x,y)\in(0,0.5]\times[-0.5,0).
	\end{cases}
\end{equation}
Solving over the domain $[-0.5,0.5]\times[-0.5,0.5]$ and assuming Dirichlet boundary conditions, this initial condition generates two shocks at angles $\theta=0$ and $\theta=\pi/2$, and two rarefaction waves at angles $\theta=\pi$ and $\theta=3\pi/2$. The solution can be found in \cite{yoon2008two}. Figure \ref{fig: burgers2d_2shocks2rare} plots the solution at final time $T_f=0.1$ using mesh $N_x=N_y=100$, and time-stepping size $\Delta t = 1.3\Delta x$ ($\text{CFL}=10.4$). The first-order scheme is compared against the third-order scheme using ENO3 and RK3.

As seen in Figure \ref{fig: burgers2d_2shocks2rare}, the $\theta_{\pi/2}$ shock collides with the $\theta_{\pi}$ rarefaction wave, and the $\theta_0$ shock collides with the $\theta_{3\pi/2}$ rarefaction wave. In the regions where the shocks interact with the rarefaction waves, there two sets of Riemann data at play: one for the shock, and another for the rarefaction wave. As with the 4 shock example, the first-order scheme slightly distorts solution around the shock-rarefaction wave interactions. Whereas, the third-order scheme better captures the sharp discontinuities without any distortion.

\begin{figure}[t!]
\centering
\begin{minipage}[b]{0.48\linewidth}
	\centering
	\includegraphics[width=\textwidth]{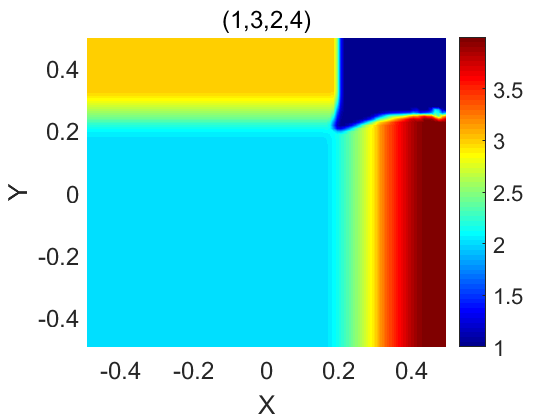}
\end{minipage}
\begin{minipage}[b]{0.48\linewidth}
	\centering
	\includegraphics[width=\textwidth]{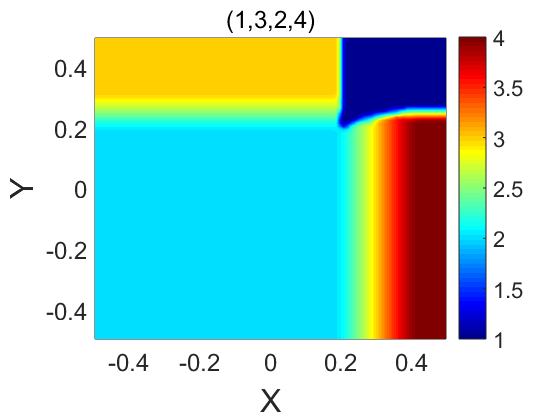}
\end{minipage}
\caption{Equation \eqref{eq: burgers2d} with initial condition \eqref{eq: 2DRiemann_2shocks2rare}: 2 shocks and 2 rarefaction waves. Mesh $N_x=N_y=100$, time-stepping size $\Delta t = 1.3\Delta x$ ($\text{CFL}=10.4$), final time $T_f=0.1$. (Left) first-order scheme, (right) third-order scheme with ENO3 and RK3.}
\label{fig: burgers2d_2shocks2rare}
\end{figure}
\end{exa}

\begin{exa}
(2D Burgers' equation with discontinuous initial data: 3 shocks and 1 rarefaction wave)\\
\ \\
We test the Riemann problem with Burgers' equation \eqref{eq: burgers2d} using initial condition \cite{yoon2008two}
\begin{equation}\label{eq: 2DRiemann_3shock1rare}
	u_0(x,y) = 
	\begin{cases}
		1,&(x,y)\in(0,0.5]\times(0,0.5],\\
		2,&(x,y)\in[-0.5,0]\times[0,0.5],\\
		3,&(x,y)\in[-0.5,0)\times[-0.5,0],\\
		4,&(x,y)\in(0,0.5]\times[-0.5,0).
	\end{cases}
\end{equation}
Solving over the domain $[-0.5,0.5]\times[-0.5,0.5]$ and assuming Dirichlet boundary conditions, this initial condition generates three shocks at angles $\theta=0$, $\theta=\pi/2$, and $\theta=\pi$, and one rarefaction wave at angle $\theta=3\pi/2$. The solution can be found in \cite{yoon2008two}. Figure \ref{fig: burgers2d_3shock1rare} plots the solution at final time $T_f=0.1$ using mesh $N_x=N_y=100$, and time-stepping size $\Delta t = 1.3\Delta x$ ($\text{CFL}=10.4$). The first-order scheme is compared against the third-order scheme using ENO3 and RK3.

As seen in Figure \ref{fig: burgers2d_3shock1rare}, the $\theta_{\pi/2}$ shock interacts with the $\theta_{\pi}$ shock at the bottom-right corner of the light blue top-left rectangle. A new shock forms at this shock-shock interaction point. Also, the $\theta_0$ shock interacts with the $\theta_{3\pi/2}$ rarefaction wave. Observe that the third-order scheme produces better resolution than the first-order scheme since it does not severely distort the solution under large time-stepping sizes.

\begin{figure}[t!]
\centering
\begin{minipage}[b]{0.48\linewidth}
	\centering
	\includegraphics[width=\textwidth]{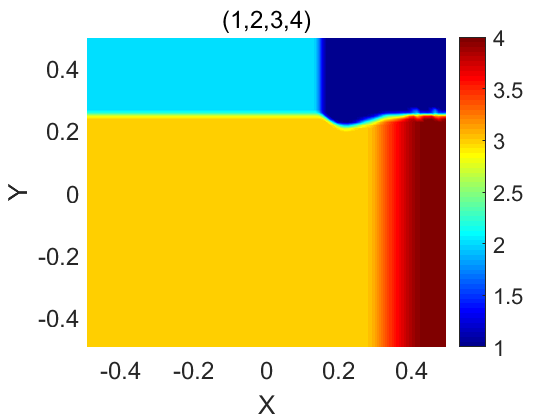}
\end{minipage}
\begin{minipage}[b]{0.48\linewidth}
	\centering
	\includegraphics[width=\textwidth]{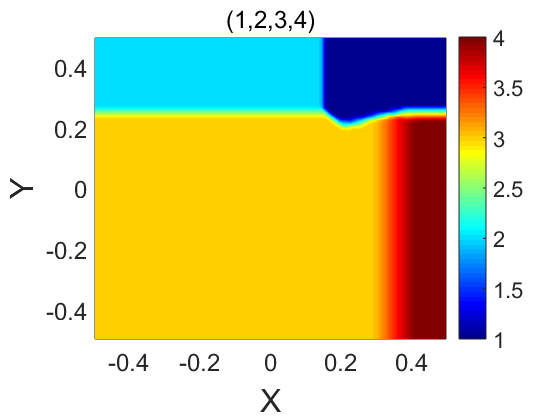}
\end{minipage}
\caption{Equation \eqref{eq: burgers2d} with initial condition \eqref{eq: 2DRiemann_3shock1rare}: 3 shocks and 1 rarefaction wave. Mesh $N_x=N_y=100$, time-stepping size $\Delta t = 1.3\Delta x$ ($\text{CFL}=10.4$), final time $T_f=0.1$. (Left) first-order scheme, (right) third-order scheme with ENO3 and RK3.}
\label{fig: burgers2d_3shock1rare}
\end{figure}
\end{exa}

\section{Conclusion}
\label{sec:conclusion}

We proposed a new high-order accurate Eulerian-Lagrangian Runge-Kutta finite volume (EL-RK-FV) scheme for solving Burgers' equation $\big(f(u) = u^2/2\big)$, specifically after shock formation. Whereas other semi-Lagrangian and Eulerian-Lagrangian schemes struggle to handle problems with intersecting characteristics, the proposed scheme remedies this issue by utilizing a merging procedure in conjuction with forward-tracing of the characteristics. ENO and WENO schemes of arbitrarily high order can be used for nonuniform spatial reconstruction, and high-order strong stability-preserving Runge-Kutta schemes can be used for the time discretization. Strang splitting was used for high-dimensional problems, but higher order splitting methods can also be used. A rigorous stability analysis for the first-order scheme was presented in \cite{yang2023stability}, and the numerical results in the current work strongly suggest that the same stability properties might hold for the high-order extension. In particular, the numerical results demonstrate our scheme's ability to achieve high-order accuracy when solving hyperbolic conservation laws after shock formation under larger time-stepping sizes compared with those in the Eulerian approach. The proposed method lays a framework that can be generalized to other nonlinear scalar conservation laws. This framework can be further combined with strategies in handling general nonlinear hypberolic systems.
\section*{Acknowledgments}
\label{sec:acknowledgments}

Research for Jing-Mei Qiu is supported by the National Science Foundation NSF-DMS-2111253, Air Force Office of Scientific Research FA9550-22-1-0390, and the Department of Energy DE-SC0023164.
\section*{Appendix: Nonuniform reconstruction procedures}
\label{sec:appendices}

\subsection*{Third-order nonuniform WENO-AO reconstruction}
Consider the indices $\mathcal{S}_j=\{j-1,j,j+1\}$. Suppose we have the cell averages $\bar{u}_i$ on nonuniform cells $I_i$, $i\in\mathcal{S}_j$. Let $s = \frac{x-x_j}{\Delta x_j}$, where $x$ is any arbitrary point, $x_j$ is the midpoint of $I_j$, and $\Delta x_j$ is the interval length associated with $I_j$. The goal is to obtain the reconstruction polynomial $P(s)$ defined over the cell $I_j$. Set linear weights as

\beq
\gamma = 0.9, \quad \gamma_l=\gamma_r=(1-\gamma)/2.
\eeq
\ \\
The reconstruction polynomial corresponding to the center stencil $\{I_{j-1},I_j,I_{j+1}\} $ is

\beq
p_c(s) = a_c \cdot s^2 + b_c\cdot s+ c_c,
\eeq
where

\beq \nonumber
\begin{split}
a_c = & \left( 3\Delta x_j^3 \Delta x_{j+1} (\Delta x_{j}+\Delta x_{j+1} )\bar{u}_{j-1} - 3 \Delta x_{j-1} \Delta x_{j}^2 \Delta x_{j+1} (\Delta x_{j-1}+2  \Delta x_{j} + \Delta x_{j+1}) \bar{u}_j  \right. \\
&\left. + 3 \Delta x_{j-1} \Delta x_{j+1}^3 (\Delta x_{j-1}+\Delta x_{j}) \bar{u}_{j+1} \right) \\
& / \left( \Delta x_{j-1} (\Delta x_{j-1}+\Delta x_{j}) \Delta x_{j+1} (\Delta x_{j}+\Delta x_{j+1}) (\Delta x_{j-1}+\Delta x_{j}+\Delta x_{j+1}) \right),
\end{split}
\eeq

\beq \nonumber
\begin{split}
b_c = & \left(-\Delta x_{j}^2 \Delta x_{j+1} (\Delta x_{j}^2 + 3 \Delta x_{j} \Delta x_{j+1} + 2 \Delta x_{j+1}^2)  \bar{u}_{j-1} \right. \\
& - \Delta x_{j-1} \Delta x_{j}  (\Delta x_{j-1}- \Delta x_{j+1} )  \Delta x_{j+1} ( 2 \Delta x_{j-1}  +  3  \Delta x_{j} + 2 \Delta x_{j+1} )\bar{u}_j \\
& \left. + \Delta x_{j-1} \Delta x_{j}^2  (2 \Delta x_{j-1}^2 + 3 \Delta x_{j-1}  \Delta x_{j} +  \Delta x_{j}^2) \bar{u}_{j+1} \right)\\
& / (\Delta x_{j-1}(\Delta x_{j-1} +\Delta x_{j}) \Delta x_{j+1} (\Delta x_{j}+\Delta x_{j+1}) (\Delta x_{j-1}+\Delta x_{j}+\Delta x_{j+1})),
\end{split}
\eeq

\beq \nonumber
\begin{split}
c_c = &\left( -\Delta x_{j}^3  \Delta x_{j+1} (\Delta x_{j}+\Delta x_{j+1}) \bar{u}_{j-1} \right.  +\Delta x_{j-1} \Delta x_{j+1} (4 \Delta x_{j-1}^2   (\Delta x_{j} + \Delta x_{j+1} ) + \\
& \Delta x_{j-1}   (3 \Delta x_{j} + 2 \Delta x_{j+1} )^2 \\
 & + \Delta x_{j} (6 \Delta x_{j}^2 + 9 \Delta x_{j} \Delta x_{j+1} + 4 \Delta x_{j+1}^2) ) \bar{u}_j \\
 &\left.- \Delta x_{j-1}  \Delta x_{j}^3 (\Delta x_{j-1}+\Delta x_{j}) \bar{u}_{j+1} \right) \\
 & / (4 \Delta x_{j-1}(\Delta x_{j-1}+\Delta x_{j}) \Delta x_{j+1} (\Delta x_{j} + \Delta x_{j+1}) (\Delta x_{j-1} + \Delta x_{j} + \Delta x_{j+1})).	
\end{split}
\eeq
\ \\
The reconstruction polynomial corresponding to the left-biased stencil $\{I_{j-1},I_j\} $ is

\beq
p_l(s) = a_l\cdot s+b_l,
\eeq
where

\beq \nonumber
a_l = \frac{2(-\Delta x_{j} \bar{u}_{j-1} \Delta x_{j}+\Delta x_{j-1} \bar{u}_{j} \Delta x_{j})} {(\Delta x_{j-1} (\Delta x_{j-1} +\Delta x_{j}))}, \quad b_l = \bar{u}_{j}.
\eeq
\ \\
The reconstruction polynomial corresponding to the right-biased stencil $\{I_j,I_{j+1}\} $ is

\beq
p_r(s) = a_r\cdot s+b_r,
\eeq
where

\beq \nonumber
a_r = \frac{2(-\Delta x_{j+1} \bar{u}_{j-1} \Delta x_{j}+\Delta x_{j} \bar{u}_{j} \Delta x_{j})} {(\Delta x_{j+1} (\Delta x_{j+1} +\Delta x_{j}))}, \quad b_r = \bar{u}_{j}.
\eeq
\ \\
The smoothness indicators are given by

\beq\nonumber
\begin{split}
   \beta_c = & \big(39 \Delta x_{j}^4 (-\Delta x_{j-1} \Delta x_{j+1} (\Delta x_{j-1} + \Delta x_{j+1}) \bar{u}_j +\Delta x_{j}^2 (\Delta x_{j+1} \bar{u}_{j-1} \\
   & + \Delta x_{j-1} \bar{u}_{j+1})+\Delta x_{j} (\Delta x_{j+1}^2 \bar{u}_{j-1} - 2 \Delta x_{j-1} \Delta x_{j+1} \bar{u}_j + \Delta x_{j-1}^2 \bar{u}_{j+1}))^2 \\
   & + \Delta x_{j+1}^2 (2 \Delta x_{j-1} \Delta x_{j+1} (\Delta x_{j-1}^2 - \Delta x_{j+1}^2 ) \bar{u}_j  + \Delta x_{j}^3 (\Delta x_{j+1} \bar{u}_{j-1} \\
   & - \Delta x_{j-1} \bar{u}_{j+1} ) + 3 \Delta x_{j+1}^2 (\Delta x_{j+1}^2 \bar{u}_{j-1} - \Delta x_{j-1}^2 \bar{u}_{j+1}) + \Delta x_{j} (2 \Delta x_{j+1}^3 \bar{u}_{j-1} \\
   & + 3 \Delta x_{j-1}^2 \Delta x_{j+1} \bar{u}_j - 3 \Delta x_{j-1} \Delta x_{j+1}^2 \bar{u}_j - 2 \Delta x_{j-1}^3 \bar{u}_{j+1}))^2 \big) \\
   & / ((\Delta x_{j-1} (\Delta x_{j-1} + \Delta x_{j}) \Delta x_{j+1} (\Delta x_{j}+\Delta x_{j+1}) (\Delta x_{j-1}+\Delta x_{j}+\Delta x_{j+1}))^2),
\end{split}
\eeq

\beq\nonumber
\beta_l = \frac{(4(-\Delta x_{j}\bar{u}_{j-1} \Delta x_{j}+\Delta x_{j-1} \bar{u}_j \Delta x_{j})^2 )} {((\Delta x_{j-1} (\Delta x_{j-1}+\Delta x_{j}))^2)},\quad 
\beta_r = \frac{(4(-\Delta x_{j+1}\bar{u}_{j} \Delta x_{j}+\Delta x_{j} \bar{u}_{j+1} \Delta x_{j})^2 )} {((\Delta x_{j+1} (\Delta x_{j+1}+\Delta x_{j}))^2)}.
\eeq
\ \\
The nonlinear weights are given by

\begin{align*}
w_c &= \gamma(1+\tau^2/(\beta_c + \epsilon)^2 ),\\
w_l &= \gamma_l (1+\tau^2/(\beta_l + \epsilon)^2 ),\\
w_r &= \gamma_r (1+\tau^2/(\beta_r + \epsilon)^2 ),
\end{align*}
where

\beq\nonumber
\tau = 0.5(\lvert \beta_c-\beta_l\rvert + \lvert \beta_c-\beta_r \rvert), \quad \text{and} \quad \epsilon = 10^{-8}.
\eeq
\ \\
Finally, the reconstructed polynomial is given by

\beq
 P(s) = \frac{w_c}{\gamma} (p_c(s) - \gamma_l \cdot p_l(s) - \gamma_r \cdot p_r(s) ) + w_l \cdot p_l(s) + w_r \cdot p_r(s).
\eeq
 
\subsection*{Third-order nonuniform ENO reconstruction}
Consider the indices $\mathcal{S}_j=\{j-2,j-1,j,j+1,j+2\}$. Suppose we have the cell averages $\bar{u}_i$ on nonuniform cells $I_i$, $i\in\mathcal{S}_j$. Let $s = \frac{x-x_j}{\Delta x_j}$, where $x$ is any arbitrary point, $x_j$ is the midpoint of $I_j$, and $\Delta x_j$ is the interval length associated with $I_j$. The goal is to choose the best 3-point stencil for the reconstruction polynomial $P(s)$ defined over the cell $I_j$,

\beq
P(s) = a\cdot s^2+ b \cdot s + c.
\eeq
\ \\
The coefficients $a$, $b$, and $c$ are chosen as follows:

\begin{itemize}
 \item If $\mathlarger{\left|\frac{\bar{u}_j -\bar{u}_{j-1}}{x_j-x_{j-1}}\right| \leq \left|\frac{\bar{u}_{j}-\bar{u}_{j+1}} {x_j-x_{j+1}} \right|}$ and $\mathlarger{\left| \frac{\mathlarger{ \frac{\bar{u}_{j-2}-\bar{u}_{j-1}}{x_{j-2}-x_{j-1}}-\frac{\bar{u}_{j-1}-\bar{u}_{j}} {(x_{j-1}-x_j)} }} {x_{j-2} - x_j} \right| \leq \left|\frac{\mathlarger{ \frac{\bar{u}_{j-1}-\bar{u}_{j}}{x_{j-1}-x_j}-\frac{\bar{u}_{j}-\bar{u}_{j+1}}{x_j-x_{j+1}} }}{(x_{j-1}-x_{j+1})} \right|}$,\\
\ \\
then we choose $\{I_{j-2}, I_{j-1}, I_j\}$ as the reconstruction stencil. The coefficients are

 \begin{subequations}
 \beq 
 \begin{split}
 a = & (3 \Delta x_{j-1} \Delta x_{j}^3  (\Delta x_{j-1}+\Delta x_{j} ) \bar{u}_{j-2} - 3 \Delta x_{j-2} \Delta x_{j}^3 (\Delta x_{j-2} + 2 \Delta x_{j-1}+ \Delta x_{j-1} ) \bar{u}_{j-1}  \\
  & + 3 \Delta x_{j-2}  \Delta x_{j-1} (\Delta x_{j-2} + \Delta x_{j-1}) \Delta x_{j}^2 \bar{u}_j) \\
  & / (\Delta x_{j-2} \Delta x_{j-1} (\Delta x_{j-2} +\Delta x_{j-1}) (\Delta x_{j-1}+ \Delta x_{j}) (\Delta x_{j-2}+\Delta x_{j-1}+\Delta x_{j}))
\end{split}
 \eeq
  \beq 
 \begin{split}
    b  = & (\Delta x_{j-1} \Delta x_{j}^2 (2 \Delta x_{j-1}^2 + 3 \Delta x_{j-1} \Delta x_{j}+\Delta x_{j}^2 ) \bar{u}_{j-2} - \Delta x_{j-2} \Delta x_{j}^2  (2  \Delta x_{j-2}^2 +6  \Delta x_{j-1}^2 + \\
      &   6 \Delta x_{j-1} \Delta x_{j} +\Delta x_{j}^2 + 3 \Delta x_{j-2} (2 \Delta x_{j-1}+\Delta x_{j})) \bar{u}_{j-1} +\Delta x_{j-2} \Delta x_{j-1} (\Delta x_{j-2} \\
      &    +\Delta x_{j-1}) \Delta x_{j} (2 \Delta x_{j-2} + 4 \Delta x_{j-1} +3 \Delta x_{j})\bar{u}_j) /(\Delta x_{j-2} \Delta x_{j-1} (\Delta x_{j-2} \\
      &    +\Delta x_{j-1}) (\Delta x_{j-1}+\Delta x_{j}) (\Delta x_{j-2} + \Delta x_{j-1} + \Delta x_{j}))
 \end{split}
 \eeq
   \beq 
 \begin{split}
    c = & (- \Delta x_{j-1} \Delta x_{j}^3 (\Delta x_{j-1} + \Delta x_{j}) \bar{u}_{j-2} +\Delta x_{j-2} \Delta x_{j}^3 (\Delta x_{j-2}+2 \Delta x_{j-1}+\Delta x_{j})\bar{u}_{j-1} \\
& +\Delta x_{j-2} \Delta x_{j-1} (4 \Delta x_{j-2}^2 (\Delta x_{j-1}+\Delta x_{j})+\Delta x_{j-1}(4 \Delta x_{j-1}^2+8 \Delta x_{j-1} \Delta x_{j}+3 \Delta x_{j}^2) \\
& +\Delta x_{j-2} (8 \Delta x_{j-1}^2 + 12 \Delta x_{j-1} \Delta x_{j}+3 \Delta x_{j}^2)) \bar{u}_j) \\
& /(4 (\Delta x_{j-2} \Delta x_{j-1} (\Delta x_{j-2}+\Delta x_{j-1}) (\Delta x_{j-1}+\Delta x_{j}) (\Delta x_{j-2}+\Delta x_{j-1}+\Delta x_{j})))
 \end{split}
 \eeq
\end{subequations}
\ \\

 \item  Else if  $\mathlarger{\left\lvert \frac{\bar{u}_j -\bar{u}_{j-1}}{x_j-x_{j-1}} \right\rvert \geq \left\lvert \frac{\bar{u}_{j}-\bar{u}_{j+1}} {x_j-x_{j+1}}\right\rvert}$ and $\mathlarger{\left\lvert \frac{\mathlarger{ \frac{\bar{u}_{j-1}-\bar{u}_{j}}{x_{j-1}-x_{j}}-\frac{\bar{u}_{j}-\bar{u}_{j+1}} {(x_{j}-x_{j+1})} }} {x_{j-1} - x_{j+1}} \right\rvert \geq \left\lvert\frac{\mathlarger{ \frac{\bar{u}_{j}-\bar{u}_{j+1}}{x_{j}-x_{j+1}} -\frac{\bar{u}_{j+1}-\bar{u}_{j+2}}{x_{j+1}-x_{j+2}}}} {(x_{j}-x_{j+2})}\right\rvert}$,\\
\ \\
then we choose $\{I_{j}, I_{j+1}, I_{j+2} \}$ as the reconstruction stencil. The coefficients are

 \begin{subequations}
\beq
\begin{split}
a = &(3 \Delta x_{j}^2 \Delta x_{j+1} \Delta x_{j+2} (\Delta x_{j+1}+\Delta x_{j+2}) \bar{u} \\ 
	  & -3 \Delta x_{j}^3 \Delta x_{j+2} (\Delta x_{j}+2 \Delta x_{j+1}+\Delta x_{j+2}) \bar{u}_{j+1} \\
	& +3 \Delta x_{j}^3 \Delta x_{j+1} (\Delta x_{j}+\Delta x_{j+1}) \bar{u}_{j+2}) \\
	& /(\Delta x_{j+1} (\Delta x_{j}+\Delta x_{j+1}) \Delta x_{j+2} (\Delta x_{j+1}+\Delta x_{j+2})(\Delta x_{j}+\Delta x_{j+1}+\Delta x_{j+2}))
\end{split}
\eeq
\beq
\begin{split}
b = & (-1 \Delta x_{j} \Delta x_{j+1} \Delta x_{j+2} (\Delta x_{j+1} + \Delta x_{j+2}) (3 \Delta x_{j}+4 \Delta x_{j+1}+2\Delta x_{j+2}) \bar{u}_j \\
& +\Delta x_{j}^2 \Delta x_{j+2} (\Delta x_{j}^2 +6 \Delta x_{j} \Delta x_{j+1}+6 \Delta x_{j+1}^2+3 \Delta x_{j}\Delta x_{j+2}+6 \Delta x_{j+1} \Delta x_{j+2} \\
& +2\Delta x_{j+2}^2) \bar{u}_{j+1} - \Delta x_{j}^2 \Delta x_{j+1} (\Delta x_{j}^2 +3 \Delta x_{j} \Delta x_{j+1}+2 \Delta x_{j+1}^2 )\bar{u}_{j+2}) \\
& /(\Delta x_{j+1} (\Delta x_{j}+\Delta x_{j+1}) \Delta x_{j+2} (\Delta x_{j+1}+\Delta x_{j+2})(\Delta x_{j}+\Delta x_{j+1}+\Delta x_{j+2}))
\end{split}
\eeq
\beq
\begin{split}
c = & (\Delta x_{j+1} \Delta x_{j+2} (\Delta x_{j+1}+\Delta x_{j+2}) (3 \Delta x_{j}^2 + 4 \Delta x_{j+1} (\Delta x_{j+1}+\Delta x_{j+2})+4 \Delta x_{j}(2 \Delta x_{j+1}\\
& +\Delta x_{j+2})) \bar{u}_j +\Delta x_{j}^3 \Delta x_{j+2} (\Delta x_{j}+2 \Delta x_{j+1}+\Delta x_{j+2}) \bar{u}_{j+1}-\Delta x_{j}^3 \Delta x_{j+1} (\Delta x_{j}+\Delta x_{j+1}) \bar{u}_{j+2}) \\
& /(4\Delta x_{j+1} (\Delta x_{j}+\Delta x_{j+1}) \Delta x_{j+2} (\Delta x_{j+1}+\Delta x_{j+2}) (\Delta x_{j}+\Delta x_{j+1}+\Delta x_{j+2}))
\end{split}
\eeq
\end{subequations}
\ \\
 \item Else, we choose $\{I_{j-1}, I_{j}, I_{j+1} \}$ as the reconstruction stencil. The coefficients are

  \begin{subequations}
 \beq 
 \begin{split}
   a = & (3 \Delta x_{j}^3 \Delta x_{j+1} (\Delta x_{j}+\Delta x_{j+1})\bar{u}_{j-1} \\
 & -3 \Delta x_{j-1} \Delta x_{j}^2 \Delta x_{j+1} (\Delta x_{j-1}+2 \Delta x_{j} +\Delta x_{j+1}) \bar{u}_j \\
 & +3 \Delta x_{j-1} \Delta x_{j}^3 (\Delta x_{j-1} + \Delta x_{j})\bar{u}_{j+1}) \\
 & /(\Delta x_{j-1} (\Delta x_{j-1}+\Delta x_{j}) \Delta x_{j+1}(\Delta x_{j}+\Delta x_{j+1}) (\Delta x_{j-1}+\Delta x_{j}+\Delta x_{j+1}))
\end{split}
 \eeq
 \beq 
\begin{split}
 b = & (-\Delta x_{j}^2 \Delta x_{j+1} (\Delta x_{j}^2 + 3 \Delta x_{j} \Delta x_{j+1} + 2 \Delta x_{j+1}^2)\bar{u}_{j-1} \\
 & -\Delta x_{j-1} \Delta x_{j} (\Delta x_{j-1} - \Delta x_{j+1}) \Delta x_{j+1} (2 \Delta x_{j-1} + 3 \Delta x_{j} + 2 \Delta x_{j+1}) \bar{u}_j  \\
 & + \Delta x_{j-1} \Delta x_{j}^2 (2 \Delta x_{j-1}^2 + 3 \Delta x_{j-1} \Delta x_{j} + \Delta x_{j}^2)\bar{u}_{j+1}) \\
 &/(\Delta x_{j-1} (\Delta x_{j-1} +\Delta x_{j}) \Delta x_{j+1} (\Delta x_{j}+\Delta x_{j+1}) (\Delta x_{j-1}+\Delta x_{j}+\Delta x_{j+1}))
\end{split}
\eeq
 \beq 
 \begin{split}
 c = & (-1 \Delta x_{j}^3 \Delta x_{j+1} (\Delta x_{j}+\Delta x_{j+1})\bar{u}_{j-1} + \Delta x_{j-1} \Delta x_{j+1} (4 \Delta x_{j-1}^2 (\Delta x_{j}+\Delta x_{j+1})+\Delta x_{j-1} \\
 & (3\Delta x_{j} +2 \Delta x_{j+1} )^2 + \Delta x_{j} (6 \Delta x_{j}^2 + 9 \Delta x_{j} \Delta x_{j+1} +4 \Delta x_{j+1}^2))\bar{u}_j \\
 & -\Delta x_{j-1} \Delta x_{j}^3 (\Delta x_{j-1}+ \Delta x_{j}) \bar{u}_{j+1}) \\
 & /(4 \Delta x_{j-1} (\Delta x_{j-1} +\Delta x_{j}) \Delta x_{j+1} (\Delta x_{j}+\Delta x_{j+1}) (\Delta x_{j-1}+\Delta x_{j}+\Delta x_{j+1}))
\end{split}
\eeq
\end{subequations}

\end{itemize}

\medskip

\printbibliography

\end{document}